\documentclass{artjlt}

\usepackage{amsmath,amsfonts,amssymb}

\title{Finite dimensional Nichols  algebras over finite cyclic groups} 
\author{Weicai Wu,  Shouchuan Zhang and  Yao-Zhong Zhang}                 
\lastname{Wu, Zhang and Zhang}  



\address{%
Weicai Wu\\               
Department  of Mathematics\\ 
Hunan University\\
Changsha  410082\\ 
P.R. China \\            
weicaiwu@hnu.edu.cn                
}
\address{%
Shouchuan Zhang\\               
Department  of Mathematics\\ 
Hunan University\\
Changsha  410082\\ 
P.R. China \\            
sczhang@hnu.edu.cn               
}
\address{%
Yao-Zhong Zhang\\               
School of Mathematics and Physics\\ 
The University of Queensland\\
Brisbane 4072\\ 
Australia \\            
yzz@maths.uq.edu.au               
}

\begin{document}


\maketitle

\begin{abstract}
All finite dimensional Nichols algebras with diagonal type of
connected finite dimensional Yetter-Drinfeld modules over finite cyclic group $\mathbb Z_n$ are found.
It is proved that Nichols algebra of connected  Yetter-Drinfeld module $V$ over $\mathbb Z_n$ 
with $\dim V >3$ is infinite dimensional.   
\end{abstract}




\numberwithin{equation}{section}

\section{Introduction}\label {s0}

This paper concerns the classification of finite dimensional pointed Hopf algebras
with finite cyclic  groups.  Recently Heckenberger established one-to-one
correspondence between arithmetic root systems and Nichols algebras of diagonal type having
a finite set of (restricted)  Poincare-Birkhoff-Witt generators \cite{He04b} and
between twisted  equivalence  classes of arithmetic root systems
and generalized Dynkin diagrams \cite{He09}. In this latter work, arithmetic root systems
were also classified in full generality.

The theory of Nichols algebras is dominated by the classification of finite dimensional
pointed Hopf algebras (see e.g. \cite {AS98, AS00}).
Nichols algebras appear in the construction of quantized Kac-Moody algebras and their
$\mathbb Z_2$-graded (see \cite {KT91,KS97}) and $\mathbb Z_3$-graded versions \cite{Ya03}.
They are natural quantum groups and are connected to the bicovariant differential calculus
initiated by Woronowicz \cite {Wo89}. Bicovariant differential calculi 
on quantum groups have been studied by Klimyk and Schm\"udgen in their book \cite{KS97}
(see especially Part IV of this book). 

Nichols algebras play a central role in the theory of (pointed) Hopf algebras. 
Any braided vector space has a canonical Nichols algebra. The
easiest braidings are those of diagonal type, that is, the vector space V
has a basis $x_1, \cdots,  x_r$ such that the braiding $c \in { \rm Aut}(V
\otimes  V )$ is given by
$c(x_i \otimes
x_j) = q_{ij}x_j \otimes
x_i $ for some nonzero numbers $q_{ij}$ , for all $i, j \in \{1,  2, \cdots,  r\}.$
The braided vector spaces of diagonal type with finite-dimensional Nichols
algebra were   essentially   classified  by Heckenberger. 
In this paper we  study diagonal braidings and their Nichols algebras coming from Yetter-
Drinfeld modules over finite cyclic groups. This is a substantial restriction,
and it turns out. We classify finite dimensional Nichols algebras with diagonal type of connected
finite dimensional Yetter-Drinfeld (YD) modules over finite cyclic group $\mathbb Z_n$.
We first determine  which braided vector space $V$ 
is a $\mathbb Z_n$-YD module by means of equation systems in $\mathbb Z_n$.
Using the classification of arithmetic root systems,
we find all finite dimensional Nichols algebras with diagonal type of connected finite dimensional $\mathbb Z_n$-YD modules.

This paper is organized as follows. In  sections 1 and 2 we find all finite dimensional
Nichols algebras with diagonal type of connected 2-dimensional and 3-dimensional
$\mathbb Z_n$-YD modules, respectively.
In section 3 we prove that Nichols algebra of connected $\mathbb Z_n$-YD module $V$ with $\dim V >3$ is infinite dimensional.

Throughout, $k$ is a  field of characteristic zero,  which contains a primitive $n$th root of unit.
 Let $G$ be a finite abelian group. Let $\widehat G := \{\chi \mid
\chi \hbox { is a homomorphism from } G \hbox { to } k^*\}$ and $R_n := \{ \omega \in  k
\mid \omega$ is a primitive $n$th root of unit$\}$.
If $G = (g)$ is a cyclic group with order $n$ and $V \in ^{kG}_{kG}{\mathcal YD }$
with basis $v_1,  v_2,  \cdots,  v_r$,  then there exist $\chi _i \in \widehat G,  g_i \in G$,
 such that $\delta  (v_i) = g_i \otimes v_i$ and $h\cdot v_i = \chi _i (h) v_i$
for  any $h \in G$,   $1\le i\le r.$
 Let $\chi\in \widehat G$ such that $\chi (g) \in R_n$.
Thus $\chi _i = \chi ^{n_i}$ and $g_i = g^{m_i}$ for $1\le i \le r.$

If $V$ is a vector space with a basis $x_1,  x_2,  \cdots,  x_r$  and  $q_{ij} \in k^*$
for $1\le i,  j \le r$ such that map $c: \left \{
\begin
{array} {lll} V \otimes V& \rightarrow   & V\otimes V \\
 x_i \otimes x_j & \mapsto   & q_{ij} x_j \otimes x_i
\end {array} \right., $ then $(V,  c)$ is called a braided vector space of diagonal type.
Denote by $(q_{ij})_{r\times r}$ the braiding matrix of $(V,  c)$ under the
basis $x_1,  x_2,  \cdots,  x_r$. Then $(V, c)$ is also written as $(V, (q_{ij})_{r \times r}).$
Let $1,  2,  \cdots,  r$ be vertexes of a diagram. There is a line between vertexes $i$ and $j$
if $q_{ij}q_{ji} \not=1$. Label vertex $i$ by $q_{ii}$ and line between $i$ and $j$ by
$q_{ij}q_{ji}$. This diagram is called  generalized Dynkin diagram (written as GDD in short)
of matrix $(q_{ij})_{r\times r}$ or $V.$  $V$ is said to be connected if the generalized
Dynkin diagram is connected. Let $e_1 := (1, 0, \cdots, 0), e_2 := (0, 1, \cdots, 0), \cdots,
e_r := (0, 0, \cdots, 1)$ be a basis of $\mathbb Z_r$. Let
$E_0 := \{e_1, e_2, \cdots, e_r\}$ and $\chi _0 (e_i, e_j) := q_{ij}.$
Then $V$ is a $\mathbb Z_r$ graded vector space if one defines ${\rm deg }\ x_i = e_i$.
Let  $$\Delta^+(\mathfrak B(V)) :=
\{ {\rm deg }\ u \mid u  \hbox { is a generator  of (restricted) PBW basis }   \}$$
and $\Delta (\mathfrak B(V)) := \Delta ^+(\mathfrak B(V))\;\cup\;-\Delta^+ (\mathfrak B(V))$.

\section{Rank 2 Nichols algebras of diagonal type}\label {s1}

In this section we find all finite dimensional Nichols algebras with diagonal type of
connected $2$-dimensional  $\mathbb Z_n$-{\rm YD} modules.

\begin{Lemma}\label{1.0} ({\rm i}) (See \cite [Lemma 2.3]{ZZC04} or appendix) Every $kG$-{\rm YD} module is a braided vector space of diagonal type.

{(\rm ii )} $V$ is a $G$-{\rm YD} module of diagonal type and braiding matrix $(q_{ij})_{n\times n}$
if and only if there exist $ \chi _j \in \widehat G$,  $g_i \in G$
such that $\chi _j (g_i) = q_{ij}$ for $1 \le i,  j \le n $.

{(\rm iii)}   If $\omega \in R_n$, then $V$ is a $\mathbb Z_n$-{\rm YD} module of diagonal type and 
 braiding matrix $(q_{ij})_{n\times n}$ if and only if there exist $m_i,  n_j \in \mathbb Z$
such that $q_{ij} = \omega ^{m_in_j}$ for $1 \le i,  j \le n $.

{(\rm iv)} If $\xi \in R_n$ and $q \in R_m$ with $m\mid n$, then there exists $s\in \mathbb Z$ such that $q = \xi ^{\frac {ns} {m}}$ with $(s, m)=1.$

{(\rm v)} If $q \in R_m$ with $m\mid n$, then there exists $\omega \in R_n$ such that $q = \omega ^{\frac {n}{m}}$.
 \end{Lemma}
\begin{Proof} ({\rm ii}) It follows from \cite [Pro. 2.4] {ZZC04}.

 ({\rm iii}) Let $G = (g)$ be a cyclic group with $\mid G \mid =n$ and $\chi \in \widehat G$ 
such that $\chi (g) =\omega$,  then $\widehat G = \{ \chi ^m \mid 1\le m \le n\}$ and 
$G = \{ g ^m \mid 1\le m \le n\}$. If $V$ is a $G$-{\rm YD} module with diagonal type 
and braiding matrix $(q_{ij})_{n\times n}$, then there exist $\chi _j \in \widehat G$ 
and $g_i \in G$ such that $\chi _j (g_i) = q_{ij}$ for $1\le i, j \le n.$ Furthermore, 
there exist $m_i, n_i$ such that $\chi _i = \chi ^{n_i}$ and $g_i = g^{m_i}$  for for $1\le i, j \le n.$  Conversely, it is clear.

 ({\rm iv}) There exist $1\le t \le n$ such that $\xi ^t = q$ with
 $m = \frac {n} {(t, n)}$. Consequently, $(t, n) = \frac {n} {m}$. There exists $s\in \mathbb Z$ such that 
$t= \frac {n}{m} s$. Let  $(s, m)=d$, $m = m'd$ and $s =s'd $. Thus $t = \frac {n}{m'} s'$.  
${\rm ord (\xi ^t)} \le m'$ since $n \mid tm'$, which implies
that $m=m'$ and $(s, m)=1.$

({\rm v}) Set $\tau := \prod \{ p \mid p $ is prime with  $p \mid n$ and $p \nmid s$ $\}$. It is clear  
$m \mid \tau$ and $(\tau +s, n)=1$. Set $\mu = \tau +s$ and $\omega := \xi ^\mu$.  Thus, $\omega ^ {\frac {n}{m}}= \xi ^ {\frac {sn}{m}}= q.$
\end{Proof}

\begin{Lemma}\label{1.1} Let $n=km$,   $(s,  m)=1,$ $t_1, t_2, t_3\in \mathbb Z.$

{(\rm i)} \begin {eqnarray}\label {e1.1.11} \left \{
\begin
{array} {lll}  x_1 y_1 & \equiv  & t_1 sk \ \ \ (\mbox {mod } n)\\
 x_2 y_2 & \equiv  & t_2 sk \ \ \ (\mbox {mod } n)\\
 x_1 y_2+ x_2y_1 & \equiv  & t_3 sk \ \ \ (\mbox {mod } n)
\end {array} \right.  \end {eqnarray}  has a solution in $\mathbb Z$ if and only if
 \begin {eqnarray}\label {e1.1.11'} \left \{
\begin
{array} {lll} \frac { x_1 y_1}{ks} & \equiv  & t_1  \ \ \ (\mbox {mod } m)\\
\frac { x_2 y_2 }{ks}& \equiv  & t_2  \ \ \ (\mbox {mod } m)\\
\frac { x_1 y_2+ x_2y_1}{ks} & \equiv  & t_3  \ \ \ (\mbox {mod } m)
\end {array} \right.  \end {eqnarray}    has a solution in $\mathbb Z$ .

{(\rm ii)}
If $d$ is  a solution of
\begin {eqnarray}\label {e1.2.1} t_1x^2-t_3x +t_2\equiv 0 \ \ \ (\mbox {mod } m),   \end {eqnarray} then
$x _1= 1,  y_1 = t_1sk,  x_2 = d,  y_2 =  (t_3-dt_1)sk$ is a solution of   (\ref  {e1.1.11})
.

({\rm iii})
If $d$ is  a solution of
\begin {eqnarray}\label {e1.2.1'} t_2x^2-t_3x +t_1\equiv 0 \ \ \ (\mbox {mod } m),   \end {eqnarray} then
 $x _2= 1,  y_2 = t_2sk,  x_1 = d,  y_1 =  (t_3-dt_2)sk$ a solution of   (\ref  {e1.1.11}).

 \end{Lemma}

\begin{Proof} It is clear. \end{Proof}

\begin{Lemma}\label{1.2} Let $n=km$ and $(s,  m)=1,$  $t_1, t_2, t_3\in \mathbb Z.$  If  (\ref {e1.1.11}) has a solution,  then
\begin {eqnarray}\label {e1.1.12} x^2-t_3x + t_1t_2\equiv 0 \ \ \ (\mbox {mod } m) \end {eqnarray} has a solution.
\end{Lemma}

\begin{Proof}
If (\ref {e1.1.11}) has a solution: $x_1= m_1, y_1= n_1, x_2= m_2, y_2= n_2$,  then
$$\left \{
\begin
{array} {lll}  \frac {s^{-1}m_1 n_1}{k} & \equiv  & t_1\ \ \ (\mbox {mod } m)\\
\frac {s^{-1}m_2 n_2}{k} & \equiv  & t_2\ \ \ (\mbox {mod } m)\\
\frac {s^{-1}m_1 n_2}{k}  + \frac {s^{-1}m_2 n_1 }{k} & \equiv  & t_3\ \ \ (\mbox {mod } m)
\end {array} \right.   \ \ \ \mbox {and} $$
$$\left \{
\begin
{array} {lll}  \frac {s^{-1}m_1 n_2}{k}\frac {s^{-1}m_2 n_1}{k} & \equiv  & t_1t_2\ \ \ (\mbox {mod } m)\\
\frac {s^{-1}m_1 n_2}{k} +\frac {s^{-1}m_2 n_1 }{k} & \equiv  & t_3\ \ \ (\mbox {mod } m)
\end {array}\right.$$ have solutions. Thus there exist $u,  v \in \mathbb Z, $ such that
$$\left \{
\begin
{array} {lll}  \frac {s^{-1}m_1 n_2}{k}\frac {s^{-1}m_2 n_1}{k} & =  & t_1t_2 + u m\ \ \ \\
\frac {s^{-1}m_1 n_2}{k} +\frac {s^{-1}m_2 n_1 }{k} & =  & t_3 + vm\ \ \
\end {array}\right.,  $$ which implies that rational number  $\frac {s^{-1}m_1n_2}{k}$ is a solution of integer coefficient equation
$ x^2-(t_3+vm)x + t_1t_2+um= 0$. Consequently,  $\frac {s^{-1}m_1n_2}{k} \in \mathbb Z.$
Therefore,  $\frac {s^{-1}m_1n_2}{k}$ is a solution of $ x^2-t_3x + t_1t_2\equiv 0\ \ \ \ \ (\mbox {mod } m)$. 
\end{Proof}

\begin{Lemma}\label{1.3} Let $n=km$ and $(s,  m)=1,$  $t_1, t_2, t_3\in \mathbb Z.$

({\rm i}) If  $m$ is odd and $(t_1, m)=1$,  then
 (\ref {e1.2.1}) has a solution if and only if  (\ref {e1.1.12}) has a solution.

  ({\rm ii}) If  $t_1$ is odd  and $(t_1, m)=1 $,  then
(\ref {e1.2.1}) has a solution if and only if  (\ref {e1.1.12}) has a solution.

({\rm iii}) If   $t_2$ is odd and $(t_2, m)=1$,  then
(\ref {e1.2.1'}) has a solution if and only if  (\ref {e1.1.12}) has a solution.

({\rm iv}) If    $(t_1, m)=1$,  then
(\ref {e1.1.11}) has a solution if and only if  (\ref {e1.2.1}) has a solution.
 \end{Lemma}
 
\begin{Proof} ({\rm i}) (\ref {e1.2.1}) and (\ref {e1.1.12}) are equivalent to
$(2t_1x-t_3 )^2 \equiv t_3^2-4t_1t_2\ \ \ (\mbox {mod } m)$ and $(2x-t_3 )^2 \equiv t_3^2-4t_1t_2\ \ \ (\mbox {mod } m)$,  
respectively. Consequently,  (\ref {e1.2.1}) has a solution if and only if  (\ref {e1.1.12}) has a solution

({\rm ii}) Considering Part {\rm (i)} we only need prove this for even $m$.
If $2\nmid t_3$ and $2 \nmid t_2$, then both  (\ref {e1.2.1}) and (\ref {e1.1.12}) have not any solutions. 
If $2\nmid t_3$ and $2 \mid t_2$,  then both (\ref {e1.2.1}) and (\ref {e1.1.12}) have  solutions. 
If $2\mid t_3$, then $(t_1x- \frac {t_3}{2})^2 \equiv (\frac {t_3}{2})^2-t_1t_2  \ \ \ (\mbox {mod } 2^{\alpha _1}) $ 
has a solution if and only if $(x- \frac {t_3}{2})^2 \equiv (\frac {t_3}{2})^2-t_1t_2  \ \ \ (\mbox {mod } 2^{\alpha _1})$ 
has a solution. Consequently, (\ref {e1.2.1}) has a solution if  and only if (\ref {e1.1.12}) has a solution.

({\rm iii}) It is similar to {\rm (ii)}.

({\rm iv}) If (\ref {e1.1.11}) has a solution, then (\ref {e1.1.12}) has a solution by Lemma \ref {1.2}, 
which implies that  (\ref {e1.2.1}) has a solution by Part ({\rm ii}).
Conversely, it follows from Lemma \ref {1.1}.  
\end{Proof}

\begin{Remark} Lemma \ref {1.2} and \ref {1.3}  hold when $s=1.$ \end{Remark}

\begin{Lemma}\label{1.4} If $(V, (q_{ij})_{r \times r})$ is a {\rm YD}- module over $\mathbb Z_n$ and  
$(V', (q_{ij}')_{r \times r})$ has degree $s_i (E_0)$ with respect to $V$ ( defined in \cite [Definition 2]{He05b}),  
then  $(V', (q_{ij}')_{r \times r})$ is also a {\rm YD}- module over $\mathbb Z_n$.
\end{Lemma}

\begin{Proof} By Lemma \ref {1.0},  there exist $m_j,  n_l \in \mathbb N$ such that $q_{jl} = \omega ^{m_jn_l}$ for $1\le j,  l\le r $,  
By \cite [Definition 2]{He05b},
  \begin {eqnarray*}q_{jl }' &=& q_{jl}q_{il} ^{m_{ij}} q_{ji} ^{m_{il}}    q_{ii} ^{m_{ij} m_{il}}  \\
 &=& \omega ^ {m_jn_l}\omega ^ {m_in_lm_{ij}} \omega ^ {m_jn_im_{il}}   \omega ^ {m_in_i m_{ij} m_{il}}\\
  &=& \omega ^ {( m_j +m_{ij}m_i )(n_l + m_{il}n_i)}.\\
   \end {eqnarray*}
Set  $m _j ':=  m_j +m_{ij}m_i,  n_l ': =n_l + m_{il}n_i.$ One has $q'_{jl} = \omega ^{m_j'n_l'}$ for $1\le j,  l \le r.$ 
Therefore,  $V'$ is  a $\mathbb Z_n$-{\rm YD} module.
\end{Proof}

Thus, if $(V, (q_{ij})_{r\times r})$ and $(V', (q_{ij}')_{r\times r})$ are Weyl equivalent, 
then $(V, (q_{ij})_{r\times r})$ is a $\mathbb Z_n$- {\rm YD} module if and only if  
$(V', (q_{ij}')_{r\times r})$ a $\mathbb Z_n$- {\rm YD} module.

\begin{Theorem}\label{1.7}  Let $n=km$ ($m>1$) and  $m = 2^{\alpha _1}3^{\alpha _2}p_3^{\alpha _3}\cdots p_r^{\alpha _r}$ 
be the prime decomposition of $m$, $r \in \mathbb N;$ $\alpha _3, \alpha _4, \cdots, \alpha _r >0,$ when $r>2$.  
If $(V, (q_{ij})_{2\times 2})$ is a braided vector space, then $V$ is a connected $\mathbb Z_n$-{\rm YD} module 
such that $\dim \mathfrak B(V)< \infty$  if and only if one of  the following conditions holds:

 T2(1).~~ $1-q_{11}q_{12}q_{21} =1-q_{12}q_{21}q_{22}=0$,  $q_{12}q_{21}\in R_m$,
  $\alpha _1=0$; $\alpha _2 =0,  1; $ $(\frac {-3}{p_i})=1$ for $ 2 < i \le r.$
Here symbol $(\frac {-3}{p_i})$ is defined in  Appendix  (\ref {e4.0.1}).

 T2(2)$_1$.~~ $1+q_{11} =1-q_{12}q_{21}q_{22}=0$,   $q_{12}q_{21}\in R_m$,
 $\alpha_1 =0$; $\alpha _1 >1$.

  T2(2)$_2$.~~  $1+q_{22} =1-q_{12}q_{21}q_{11}=0$, $q_{12}q_{21}\in R_m$,
 $\alpha_1 =0$; $\alpha _1 >1$.

  T2(3).~~  $1+q_{11} =1+q_{22}=0$,  $q_{12}q_{21}\in R_m$,
   $\alpha_1 =0  $; $\alpha _1 >1$.

  T3(1)$_1$.~~   $q_{12}q_{21} =q_{11}^{-2}$,  $q_{22} = q_{11}^2$,  $q_{11} \in R_m$,  $m >2$; 
$\alpha _1=0,  1;$  $\alpha _2=0$; $p_i \equiv 1 \ \ \ (\mbox {mod } 4), $  for $ 2 < i \le r.$

 T3(1)$_2$.~~ $q_{12}q_{21} =q_{11}^{-2}$,  $q_{22} =-1$, $ q_{11}\in {R_{m}}$, $m>2$,
 $\alpha _1 \not= 2,  3$.

 T3(2)$_1$.~~  $\omega \in R_n$, $s =1, 2;$ $q_{11}= \omega ^{\frac {ns}{3}}$,
 $q_{22}= \omega ^{\frac {n}{m}}$, $q_{12}q_{21} q_{22}=1$, $m>3;$ $3\nmid m$ or $ \frac {ms}{3} \not\equiv 2  \ \ \ (\mbox {mod } 3)$.

T3(2)$_2$.~~  $q_{12}q_{21} q_{22}=1$,  $q_{11}\in{R_{3}}$,
$q_{22}\in{R_{2}}$,
   $m=6 $;

 T3(3).~~ $q_{11}\in{R_{3}}$, $q_{12}q_{21}=-q_{11}$, $q_{22}=-1$;
 $m=6.$

  T4(1).~~ $q_{0}=q_{12}q_{21}q_{11}\in R_{12}$, $q_{11} = q_{0}^4$, $q_{22} =-q_{0}^2$,
  $m = 12$.

T4(2).~~ $q_{12}q_{21}\in R_{12}$, $q_{11} =
q_{22}=-(q_{12}q_{21})^2$, $m = 12$.

T5(1).~~  $q_{12}q_{21}\in R_{12}$, $q_{11} =-(q_{12}q_{21})^2$, $q_{22}=-1$, $m = 12$.

T5(2).~~  $q_{0}=q_{12}q_{21}q_{11}\in R_{12}$, $q_{11} =
q_{0}^4$, $q_{22}=-1$, $m = 12$.

T6.~~  $q_{11}\in{R_{18}}$, $ q_{12}q_{21}=q_{11}^{-2}$, $
q_{22}=-q_{11}^3$,  $m = 18$.

T7(1).~~ $q_{11}\in R_{12}$, $ q_{12}q_{21}=q_{11}^{-3}$, $q_{22}=-1$; $m = 12$.

T7(2).~~ $q_{12}q_{21}\in R_{12}$, $q_{11}
=(q_{12}q_{21})^{-3}$, $q_{22}=-1$; $m = 12$.

 T8(1).~~ $ q_{12}q_{21}=q_{11}^{-3}$,  $
q_{22}=q_{11}^3$, $q_{11} \in R_{m}$, $m>3$,  $\alpha _1=0$; $\alpha _2 =0,  1; $ $(\frac {-3}{p_i})=1$ for  $2 < i \le r$.

 T8(2).~~$_{1}$ $ (q_{12}q_{21})^4=-1$, $q_{22}=-1$, $q_{12}q_{21}=-q_{11}$;  $m=8$.

  T8(2).~~$_{2}$ $ (q_{12}q_{21})^4=-1$, $q_{22}=-1$, $q_{11} =(q_{12}q_{21})^{-2}$;  $m = 8$.

 T8(3).~~ $ (q_{12}q_{21})^4=-1$, $q_{11} =(q_{12}q_{21})^2$, $q_{22} =(q_{12}q_{21})^{-1}$; $m = 8$.

 T9.~~  $q_{12}q_{21}\in R_{9}$, $q_{11}
=(q_{12}q_{21})^{-3}$, $q_{22}=-1$; $m = 18$.

 T10.~~  $q_{12}q_{21}\in R_{24}$, $q_{11} =(q_{12}q_{21})^{-6}$, $q_{22} =(q_{12}q_{21})^{-8}$; $m = 24$.

 T11(1).~~  $q_{11}\in R_{5}$, $ q_{12}q_{21}=q_{11}^{-3}$, $q_{22}=-1$; $m =10$.

T11(2).~~  $q_{11}\in R_{20}$, $q_{12}q_{21}=q_{11}^{-3}$, $q_{22}=-1$; $m= 20$.

 T12.~~   $q_{11}\in{R_{30}}$, $ q_{12}q_{21}=q_{11}^{-3}$, $
q_{22}=-q_{11}^5$; $m=30$.

 T13.~~  $q_{12}q_{21}\in R_{24}$, $q_{11} =(q_{12}q_{21})^6$, $q_{22} =(q_{12}q_{21})^{-1}$; $m= 24$.

 T14.~~   $q_{11}\in R_{18}$, $ q_{12}q_{21}=q_{11}^{-4}$, $q_{22}=-1$; $m=18$.

 T15.~~  $q_{12}q_{21}\in R_{30}$, $q_{11} =-(q_{12}q_{21})^{-3}$, $q_{22} =(q_{12}q_{21})^{-1}$; $m=30$.

T16(1).~~ $q_{11}\in R_{10}$, $ q_{12}q_{21}=q_{11}^{-4}$, $q_{22}=-1$; $m=10$.

 T16(2).~~ $q_{12}q_{21}\in R_{20}$, $q_{11} =(q_{12}q_{21})^{-4}$, $q_{22}=-1$; $m=20$.

 T17.~~  $q_{12}q_{21}\in R_{24}$, $q_{11} =-(q_{12}q_{21})^4$, $q_{22}=-1$; $m=24$.

 T18.~~  $q_{12}q_{21}\in R_{30}$, $q_{11} =-(q_{12}q_{21})^5$, $q_{22}=-1$; $m=30$.

 T20.~~  $q_{12}q_{21}\in R_{30}$, $q_{11} =(q_{12}q_{21})^{-6}$, $q_{22}=-1$; $m=30$.

 T21.~~  $q_{11}\in R_{24}$, $ q_{12}q_{21}=q_{11}^{-5}$, $q_{22}=-1$; $m= 24$.

 \end{Theorem}

\begin{Proof}  By \cite [Th. 4] {He04a},  it is enough to check if  there exist
$\mathbb Z_n$-{\rm YD} satisfying T2-T22.

T2$(1)$
If   \begin {eqnarray}\label {e1.1.131} \left \{
\begin
{array} {lll}  x_1 y_1 & \equiv  & \frac {sn}{m} \ \ \ (\mbox {mod } n)\\
 x_2 y_2 & \equiv  & \frac {sn}{m} \ \ \ (\mbox {mod } n)\\
 x_1 y_2+ x_2y_1 & \equiv  & -\frac {sn}{m} \ \ \ (\mbox {mod } n)
\end {array} \right.,  \end {eqnarray} has a solution, where  $(s,  m)=1$,  then
 $
x^2 +x + 1\equiv 0 \ \ \ (\mbox {mod } m) $ has a solution,  which implies $\alpha _1 =0$ and
$ (2x+ 1)^2 \equiv-3 \ \ \ (\mbox {mod } p_i^{\alpha _i} ) $ has a solution for
$2< i \le r$. It is clear that $ (2x+ 1)^2 \equiv-3 \ \ \ (\mbox {mod } 3 ) $ has a solution and 
$ (2x+ 1)^2 \equiv-3 \ \ \ (\mbox {mod } 3^2 ) $ has not any solution.
thus   $\alpha _2 =0,  1; $ $(\frac {-3}{p_i})=1$ for
$2< i \le r$. Conversely, (\ref {e1.1.11}) has a solution by Lemma \ref {1.1} since (\ref {e1.2.1}) 
has a solution when  $\alpha _1=0$; $\alpha _2 =0,  1; $ $(\frac {-3}{p_i})=1$ for
$2< i \le r$.

 T2(2)$_1$ ({\rm i}) $2\mid m.$  If  \begin {eqnarray} \label {e2.6.3} \left \{
\begin
{array} {lll}  x_1 y_1 & \equiv  & \frac {sn}{m} \ \ \ (\mbox {mod } n)\\
 x_2 y_2 & \equiv  & \frac {sn}{2} \ \ \ (\mbox {mod } n)\\
 x_1 y_2+ x_2y_1 & \equiv  & -\frac {sn}{m} \ \ \ (\mbox {mod } n)
\end {array} \right.,  \end {eqnarray} has a solution, where  $(s,  m)=1$,  then
 $
x^2 +x + \frac {m}{2}\equiv 0 \ \ \ (\mbox {mod } m) $ has a solution, which implies $\alpha _1>1$ and
$ (2x+ 1)^2 \equiv 1 \ \ \ (\mbox {mod } p_i^{\alpha _i} ) $  by Lemma \ref {4.2} ({\rm i}) for  $1 <i \le r.$

 ({\rm ii})  $2\nmid m$ and $2\mid n$.  Since
 $
mx^2 +2sx +2s\equiv 0 \ \ \ (\mbox {mod } 2m) $  has always a solution,  \begin {eqnarray*} \left \{
\begin
{array} {lll}  x_1 y_1 & \equiv  & 2sk_1 \ \ \ (\mbox {mod } n)\\
 x_2 y_2 & \equiv  & mk_1 \ \ \ (\mbox {mod } n)\\
 x_1 y_2+ x_2y_1 & \equiv  & -2sk_1 \ \ \ (\mbox {mod } n)
\end {array} \right.  \end {eqnarray*} has a solution by Lemma \ref {1.1} ({\rm ii}),  where  $(s,  m)=1$ and $n = 2mk_1$.

 T2 (2)$_2$  It is similar to T2 (2)$_1$ .

 T2 $(3)$ ({\rm i}) $2\mid m$. Considering  Lemma \ref {1.1}({\rm i}) one obtains that   \begin {eqnarray*} \left \{
\begin
{array} {lll}  x_1 y_1 & \equiv  & \frac {n}{2} \ \ \ (\mbox {mod } n)\\
 x_2 y_2 & \equiv  & \frac {n}{2} \ \ \ (\mbox {mod } n)\\
 x_1 y_2+ x_2y_1 & \equiv  &  \frac {sn}{m} \ \ \ (\mbox {mod } n)
\end {array} \right.,  \end {eqnarray*} has not any solution, where  $(s,  m)=1$,
since $ x^2-x + \frac {m^2}{4}  \equiv   0  \ \ \ (\mbox {mod } m)$ has  not any solutions when $\alpha_1=1$ 
by Lemma \ref {4.2}({\rm i}). It is clear that
$\frac {m}{2}x^2-x +\frac {m}{2}\equiv 0 \ \ \ (\mbox {mod } 2^{\alpha_1})$ has a
solution $2^{\alpha _1-1}$ when $\alpha_1 >1.$

({\rm ii})  $2\nmid m$ and $2\mid n$.  One obtains that  \begin {eqnarray*} \left \{
\begin
{array} {lll}  x_1 y_1 & \equiv  & mk_1 \ \ \ (\mbox {mod } n)\\
 x_2 y_2 & \equiv  & mk_1 \ \ \ (\mbox {mod } n)\\
 x_1 y_2+ x_2y_1 & \equiv  &  2sk_1 \ \ \ (\mbox {mod } n)
\end {array} \right.,  \end {eqnarray*} has a solution, where  $(s,  m)=1$ and $n = 2mk_1$, since
 $
mx^2-2sx +m\equiv 0 \ \ \ (\mbox {mod } 2m) $   has always a solution.

 T3 (1)$_1$  By Lemma \ref {1.1}({\rm i}) and  \begin {eqnarray*} \left \{
\begin
{array} {lll}  x_1 y_1 & \equiv  & \frac {sn}{m} \ \ \ (\mbox {mod } n)\\
 x_2 y_2 & \equiv  & \frac {2sn}{m} \ \ \ (\mbox {mod } n)\\
 x_1 y_2+ x_2y_1 & \equiv  &  \frac {-2sn}{m} \ \ \ (\mbox {mod } n)
\end {array} \right.,  \end {eqnarray*} where  $(s,  m)=1$,  one obtains
\begin {eqnarray*}
x^2 +2x +2\equiv 0 \ \ \ (\mbox {mod } m) \end {eqnarray*} and  \begin {eqnarray*}
(x+1)^2\equiv-1 \ \ \ (\mbox {mod } m) \end {eqnarray*} which implies $\alpha _1 = 0,  1; \alpha _2=0;$ $(\frac {-1}{p_i})=1$ for $2 <i \le r.$

 T3 (1)$_2$  ({\rm i}) $2 \mid m.$    By Lemma \ref {1.1}({\rm i}) and  \begin {eqnarray*} \left \{
\begin
{array} {lll}  x_1 y_1 & \equiv  & \frac {sn}{m} \ \ \ (\mbox {mod } n)\\
 x_2 y_2 & \equiv  & \frac {n}{2} \ \ \ (\mbox {mod } n)\\
 x_1 y_2+ x_2y_1 & \equiv  &  \frac {-2sn}{m} \ \ \ (\mbox {mod } n)
\end {array} \right.,  \end {eqnarray*} where  $(s,  m)=1$, one obtains
\begin {eqnarray*}
x^2 +2x +\frac {m}{2}\equiv 0 \ \ \ (\mbox {mod } m) \end {eqnarray*} and  \begin {eqnarray*}
(x+1)^2\equiv 1-\frac {m}{2} \ \ \ (\mbox {mod } m). \end {eqnarray*}
By Lemma \ref {4.3}({\rm ii}),  \begin {eqnarray*}
(x+1)^2\equiv 1-\frac {m}{2} \ \ \ (\mbox {mod } 2^2) \end {eqnarray*} and   \begin {eqnarray*}
(x+1)^2\equiv 1-\frac {m}{2} \ \ \ (\mbox {mod } 2^3)\end {eqnarray*} has not any solutions.
\begin {eqnarray*}
(x+1)^2\equiv 1-\frac {m}{2} \ \ \ (\mbox {mod } 2^{\alpha_1}) \end {eqnarray*}  has a solution when $\alpha _1>3.$

  ({\rm ii}) $2 \nmid m.$ $n = 2mk_1.$ By Lemma \ref {1.1}({\rm i}) and  \begin {eqnarray*} \left \{
\begin
{array} {lll}  x_1 y_1 & \equiv  & 2sk_1 \ \ \ (\mbox {mod } n)\\
 x_2 y_2 & \equiv  &  mk_1 \ \ \ (\mbox {mod } n)\\
 x_1 y_2+ x_2y_1 & \equiv  & -4sk_1 \ \ \ (\mbox {mod } n)
\end {array} \right.,  \end {eqnarray*} where  $(s,  m)=1$,   one obtains
\begin {eqnarray*}
mx^2 +4sx +2s\equiv 0 \ \ \ (\mbox {mod } 2m) \end {eqnarray*}  has a solution.

T3 $(2)_1$ ({\rm i}) $3\mid m$. By Lemma \ref {1.0}, one has
\begin {eqnarray}\label {e1.1.131a} \left \{
\begin
{array} {lll}  x_1 y_1 & \equiv  & \frac {sn}{3}  \ \ \ (\mbox {mod } n)\\
 x_2 y_2 & \equiv  & \frac {n}{m}  \ \ \ (\mbox {mod } n)\\
 x_1 y_2+ x_2y_1 & \equiv  &- \frac {n}{m} \ \ \ (\mbox {mod } n)
\end {array} \right..  \end {eqnarray}  Let $m = 3m'.$
 $
x^2 +x + \frac {sm}{3}\equiv 0 \ \ \ (\mbox {mod } m).$
 $
(2x+1)^2\equiv 1-4 m's \ \ \ (\mbox {mod } p_i ^{\alpha _i})$ has a solution for $
2<i \le r $.
    Consequently,  $ 1-4 m's \not\equiv 2 \ \ \ (\mbox {mod } 3 )$ since $ m's \not\equiv 2 \ \ \ (\mbox {mod } 3 )$. 
This implies $  (2x+1)^2\equiv  1-4 m's  \ \ \ (\mbox {mod } 3 )$ has a solution.

({\rm ii}) $3\nmid m$. $n=mk$ and $k =3k_1$.   If  \begin {eqnarray} \label {e1.8.22} \left \{
\begin
{array} {lll}  x_1 y_1 & \equiv  & \frac {s_1n}{3}  \ \ \ (\mbox {mod } n)\\
 x_2 y_2 & \equiv  & \frac {sn}{m}  \ \ \ (\mbox {mod } n)\\
 x_1 y_2+ x_2y_1 & \equiv  &- \frac {sn}{m} \ \ \ (\mbox {mod } n)
\end {array} \right.  \end {eqnarray} has a solution, where $s_1=1$ or $2$,  $(s,  m)=1$. 
then  $ms_1x^2 + 3sx + 3s \equiv 0 \ \ \ (\mbox {mod } 3m)$ has a solution, 
which implies that (\ref {e1.8.22}) has a solution by Lemma \ref {1.1}.

T3(2)$_2$   (\ref {e1.2.1}) has a solution $d=3$ with  $m=6,$ $t_1 = 2,  $ $t_2 = 3, $ $ t_3 =-3$.

T3(3) (\ref {e1.2.1}) has a solution $d=1$ with $m=6,$ $t_1 = 2,  $ $t_2 = 3, $ $ t_3 = 5 $.

T4(1)  (\ref {e1.2.1}) has a solution $d=4$ with $m =12$,  $t_1 = 4,  $ $t_2 = 8, $ $ t_3 = 9.$

T4(2) (\ref {e1.2.1}) has a solution $d=4$ with $m =12$,  $t_1 = 8,  $ $t_2 = 8, $ $ t_3 = 1.$

T5(1)  (\ref {e1.2.1}) has a solution $d=2$ with $m =12$,  $t_1 = 8,  $ $t_2 = 6, $ $ t_3 = 1.$

T5(2)  (\ref {e1.2.1}) has a solution $d=6$ with $m =12$,  $t_1 = 4,  $ $t_2 = 6, $ $ t_3 = 9.$

T6 (\ref {e1.2.1}) has a solution $d=12$ with $m =18$,  $t_1 = 1,  $ $t_2 = 12, $ $ t_3 = 16.$

T7(1)  (\ref {e1.2.1}) has a solution $d=3$ with $m =12$,  $t_1 = 1,  $ $t_2 = 6, $ $ t_3 =-3.$

T7(2)  (\ref {e1.2.1}) has a solution $d=3$ with $m =12$,  $t_1 =-3,  $ $t_2 = 6, $ $ t_3 = 1.$

T8 (1) By  \begin {eqnarray*} \left \{
\begin
{array} {lll}  x_1 y_1 & \equiv  & \frac {sn}{m} \ \ \ (\mbox {mod } n)\\
 x_2 y_2 & \equiv  & \frac {3sn}{m} \ \ \ (\mbox {mod } n)\\
 x_1 y_2+ x_2y_1 & \equiv  &  \frac {-3sn}{m} \ \ \ (\mbox {mod } n)
\end {array} \right.,  \end {eqnarray*} where  $(s,  m)=1$,  one obtains
\begin {eqnarray} \label {e1.6.2}
x^2 +3x +3\equiv 0 \ \ \ (\mbox {mod } m). \end {eqnarray}
By Lemma \ref {4.2},  $x^2 +3x +3\equiv 0 \ \ \ (\mbox {mod } 2)$ does not have
 any solutions, which implies $\alpha_1 =0$;
$(2x+3)^2 \equiv-3 \ \ \ (\mbox {mod } p^{\alpha_i})$ for $ 2 < i\le r$, which implies $(\frac {-3}{p_i}) =1$ 
for $ 2 < i\le r$;  $(2x+3)^2 \equiv-3 \ \ \ (\mbox {mod } 3)$ has a solution and 
$(2x+3)^2 \equiv-3 \ \ \ (\mbox {mod } 3^2)$ does not have any solutions,
 which implies $\alpha _2 =0, 1$.

T8(2)$_1$  (\ref {e1.2.1}) has a solution $d=4$ with $m =8$,  $t_1 = 5,  $ $t_2 = 4, $ $ t_3 = 1.$

T8(2)$_2$    (\ref {e1.2.1}) has a solution $d=4$ with $m =8$,  $t_1 =-2,  $ $t_2 = 4, $ $ t_3 = 1.$

T8(3)    (\ref {e1.2.1}) has a solution $d=1$ with $m =8$,  $t_1 = 2,  $ $t_2 =-1, $ $ t_3 = 1.$

T9  (\ref {e1.2.1'}) has a solution $d=6$ with $m =18$,  $t_1 =-6,  $ $t_2 = 9, $ $ t_3 = 2.$

T10  (\ref {e1.2.1}) has a solution $d=16$ with $m =24$,  $t_1 =-6,  $ $t_2 =-8, $ $ t_3 = 1.$

T11(1)    (\ref {e1.2.1'}) has a solution $d=8$ with $m =10$,  $t_1 = 2,  $ $t_2 = 5, $ $ t_3 =-6.$

T11 (2)  (\ref {e1.2.1}) has a solution $d=7$ with $m =20$,  $t_1 = 1,  $ $t_2 = 10, $ $ t_3 =-3.$

T12  (\ref {e1.2.1}) has a solution $d=10$ with $m =30$,  $t_1 = 1,  $ $t_2 = 20, $ $ t_3 =-3.$

T13  (\ref {e1.2.1}) has a solution $d=5$ with $m =24$,  $t_1 = 6,  $ $t_2 =-1, $ $ t_3 = 1.$

T14  (\ref {e1.2.1}) has a solution $d=9$ with $m =18$,  $t_1 = 1,  $ $t_2 = 9, $ $ t_3 =-4.$

T15  (\ref {e1.2.1}) has a solution $d=11$ with $m =30$,  $t_1 = 12,  $ $t_2 =-1, $ $ t_3 = 1.$

T16 (1)  (\ref {e1.2.1}) has a solution $d=5$ with $m =10$,  $t_1 = 1,  $ $t_2 = 5, $ $ t_3 =-4.$

T16 (2)  (\ref {e1.2.1}) has a solution $d=10$ with $m =20$,  $t_1 =-4,  $ $t_2 = 10, $ $ t_3 = 1.$

T17  (\ref {e1.2.1}) has a solution $d=4$ with $m =24$,  $t_1 = 16,  $ $t_2 = 12, $ $ t_3 = 1.$

T18  (\ref {e1.2.1}) has a solution $d=5$ with $m =30$,  $t_1 = 20,  $ $t_2 = 15, $ $ t_3 = 1.$

 T19  (\ref {e1.1.12}) becomes $ x^2 + 3 x +7  \equiv 0\ \ \ (\mbox {mod } 14)$,
which does not have any solution by Lemma \ref {4.2}({\rm i}).

T20  (\ref {e1.2.1}) has a solution $d=15$ with $m =30$,  $t_1 =-6,  $ $t_2 = 15, $ $ t_3 = 1.$

T21  (\ref {e1.2.1}) has a solution $d=7$ with $m =24$,  $t_1 = 1,  $ $t_2 = 12, $ $ t_3 =-5.$

T22  (\ref {e1.1.12}) becomes $ x^2 + 5 x +7 \equiv 0\ \ \ (\mbox {mod } 14)$,
which does not have any solutions by Lemma \ref {4.2}({\rm i}).
\end{Proof}

\begin{Proposition}\label{1.6}  If $(V, (q_{ij})_{2\times 2})$ is a braided vector space and $q_{ij}$ 
is a root of unit for $i, j =1, 2$, then  $\dim \mathfrak B(V)< \infty$  if and only if $\Delta (\mathfrak B(V))$ is finite.
\end {Proposition}

\begin{Proof}  It is clear that  $\Delta (\mathfrak B(V))$
is finite if $\dim \mathfrak B(V)< \infty$ by \cite {He06b}. Conversely,  if  $\Delta (\mathfrak B(V))$
is finite,  then the generalized Dynkin diagram of $V$ is in \cite [Table 1]{He05c}.
It follows  $\dim \mathfrak B(V)< \infty$ from \cite [Th. 4] {He04a}. 
\end{Proof}

\section{Rank 3 Nichols algebras of diagonal type}\label {s2}

In this section we present all finite  dimensional Nichols algebras with diagonal type of
connected $3$-dimensional $\mathbb Z_n$-{\rm YD} modules.

Let $|u|$ denote length of word $u$.
 \begin{Lemma}\label{2.1'}
 {(\rm i)} If $\mid u \mid = \mid v \mid$, then $u<v$ if and only if $uw < vw.$

 {(\rm ii)} If $u=vw$ is the Shirshow decomposition of Lyndon word $u$  and  $[u]$ is  hard,   then both  $[v]$ and $[w]$ are  hard too.
 \end{Lemma}
 
\begin{Proof}
 ({\rm i}) It is clear.

 ({\rm ii}) If  $[w]$ is not  hard,  then  there exist words $w_i> w$ and $k_i \in k$ such that 
$w =\sum _{i=1}^m k_i w_i$ by \cite [Cor. 3.2.4] {Kh99b}. Consequently,  $u= vw = \sum _{i=1}^m k_i vw_i$ 
and $[u]$ is not a hard word by \cite [Cor. 3.2.4] {Kh99b}. This is a contradiction. If  $[v]$ is not  hard,  
then  there exist words $v_i> v$ and $k_i \in k$ such that $v =\sum _{i=1}^m k_i v_i$ by \cite [Cor. 3.2.4] {Kh99b}. 
Consequently,  $u= vw = \sum _{i=1}^m k_i v_iw$ and  $v_iw > vw$ by Part ({\rm i}), 
which implies that  $[u]$ is not a hard word by \cite [Cor. 3.2.4] {Kh99b}.  
\end{Proof}

Let $\chi _u $ and $g_u$ denote $\chi _{i_1} * \chi _{i_2} * \cdots *\chi _{i_r}$  and 
$g _{i_1}  g _{i_2}  \cdots g _{i_r}$,   respectively,    for  any homogeneous
 element $u \in \mathfrak B(V) $ with  $deg (u)  = g_{i_1} g_{i_2} \cdots g_{i_r}$,   
where $(\chi _{i_1} * \chi _{i_2} * \cdots \chi _{i_r}) (g)  = \chi _{i_1} (g)  \chi _{i_2}(g)  \cdots \chi _{i_r}(g) $.
Define
\begin {eqnarray}\label{define} 
[u,   v] = vu   -p_{v,   u}uv 
\end {eqnarray}
 and $ [u,   v]_c =  [v,   u] $,   where  $p_{u,  v} = \chi _v(g_u) $.  By \cite {ZZ04},
$(\mathfrak B(V),    [\ ]_c) $ is a braided  m-Lie algebra and we have the braided Jacobi identity as follows:
\begin {eqnarray}\label {e2.1.2}   
[ [u,  v],  w]= [u,   [v,  w]]  +p_{vw}^{ -1} [ [u,  w],  v]  +(p_{wv} -p_{vw}^{ - 1}) v\cdot [u,  w].
\end {eqnarray}

Recall duality $\mathfrak B(V^*) $ of Nichols algebra $\mathfrak B(V) $ in \cite [Section 1.3]{He05} 
and \cite{He06b}. Let $y_1,   y_2,   y_3$ be a dual basis of $x_1,   x_2,   x_3$. $\delta (y_i)  = g_i ^{-1} \otimes y_i$,   
$g_i \cdot y_j = q_{ij}^{-1} y_j $ and $\Delta (y_i)  = g_i ^{-1} \otimes y_i +y_i \otimes 1.$ There exists a bilinear
 map $<, >$ from $(\mathfrak B(V^*)  \# kG)  \times \mathfrak B(V) $ to $\mathfrak B(V) $ such that     
$<y_i,   uv> = <y_i,   u>v +g_i^{-1}.u<y_i,   v>$ and
$<y_i,   <y_j,   u>> = <y_iy_j,   u>$  for any $u,   v\in  \mathfrak B(V)  $ and  $i=1, 2, 3$. 
Furthermore,   for any $u\in \oplus _{i=1}^\infty \mathfrak B(V)_{(i)}$,   one has that  $u=0$ if and only if 
$<y_i,   u> = 0$ for $i =1,   2,   3.$ We often use this to show many relations.

 Let $1,   2,   3$ denote $x_1,   x_2,   x_3$ in short,   respectively.
\begin{Lemma}\label{2.2'''} Let $q_{11} = -1$,   $q_{23}q_{32}=1$. Then

{\rm (i)}   1)  $<y_k,   [j,   i]>=0,  \forall\ k\neq j$.

2)  $ [ [1,  3],  2]=q_{32}^{ -1} [ [1,  2],  3]$,   $<y_i,   [ [1,  3],  2]>=0,  $ for $i=2,  3$.

3)  $[2, 3] =0$ and   $32=q_{32}23$.

4)  $[1,   [1,   2]] = [1,   [1,   3]]=0$.

{\rm (ii)}  $<y_1,    [  [1,  3],  2]>
=(q_{12}^{ -1} -q_{21}) (q_{13}^{ -1} -q_{31}) 23$.

{\rm (iii)}  $<y_1,   [ [1,  2],   [1,  3]]>  = -q_{12}^{ -1}q_{13}^{ -1}(1 -q_{12}q_{21}q_{31}q_{13})   [2,    [1,  3]]$

$=q_{13}^{ -1}(1 -q_{12}q_{21}q_{31}q_{13}) (q_{32}231 -q_{12}^{ -1}312  +q_{12}^{ -1}q_{31}q_{32}123 - q_{31}q_{32}213) $.

{\rm (iv)}   $ <y_1,   [ [ [1,  2],      [1,  3]],  2]>=  -q_{12}^{
-1}q_{13}^{ -1}(1 -q_{12}q_{21}q_{31}q_{13}) \left(q_{12}^{ -1}2[2  [1, 3]] \right.$  

$ \left.-q_{21}^2q_{22}q_{23}[2  [1,  3]]2\right) $.

{\rm (v)}  Furthermore,   if $(q_{22} +1) (q_{22} q_{12}q_{21}-1) = (q_{33} +1) (q_{33} q_{13}q_{31}-1) =0  $,   then

1)  $[[1,   2],   2] =[[1,   3],   3] =0$.

2)  $ [ [1,  2],   [ [1,  3],  2]]= [ [ [1,  2],   [1,  3]],  2]$.

{\rm (vi)}  Furthermore,  if  $q_{22} = q_{33} =-1$, then
 $<y_1,    [  [1,  3],    [  [1,  3],  2]]> $

$= \{ -(q_{12}^{ -1} -q_{21}) (q_{13}^{ -1} -q_{31}) q_{31}+q_{11}^{ -1}q_{13}^{ -1}q_{12}^{ -1}(q_{13}^{ -1}-q_{31})$

$+q_{11}q_{13}q_{31}q_{33}q_{21}(q_{13}^{ -1} -q_{31}) q_{31}\}2313  
+\{ -q_{11}^{ -1}q_{13}^{ -1}q_{12}^{ -1}q_{21}q_{23}(q_{13}^{ -1} -q_{31})$

$ - q_{11}q_{13}q_{31}q_{33}q_{21}q_{23}(q_{13}^{ -1}
 -q_{31}) q_{21}q_{31} -q_{31}q_{33}q_{21}q_{23}(q_{12}^{ -1} - q_{21}) (q_{13}^{ -1} -q_{31}) \}3123.$

{\rm (vii)}   Furthermore,  if  $q_{22} = q_{33} =-1$, then  $<y_1,    [  [  [1,  2],    [1,  3]],    [1,  3]]>$

$=q_{13}^{ -2}q_{12}^{ -1}(1 -q_{12}q_{21}q_{31}q_{13}-q_{12}q_{21}q_{31}q_{13}q_{33}
 +q_{12}q_{21}q_{31}^{2}q_{13}^{2}q_{33})   [1,  3]^{2}2$

$+q_{12}q_{32}^{2}q_{13}^{ -1}q_{31}(1 -q_{31}q_{13} -
q_{31}q_{13}q_{33}  +q_{12}q_{21}q_{31}^{2}q_{13}^{2}q_{33}) 2  [1, 3]^{2} $

$ +q_{32}q_{13}^{ -2}( -1  +q_{31}q_{13}q_{33}
+q_{12}q_{21}q_{31}^{2}q_{13}^{2} - q_{12}q_{21}q_{31}^{3}q_{13}^{3}q_{33})   [1,  3]2  [1,  3]$.

{\rm (viii)}    Furthermore,  if  $q_{22} = q_{33} =-1$, then  $ <y_1,   [  [1, 2],   [  [1, 2],   [1, 3]]]> $

$=\{(q_{12}^{ -1} -q_{21})  -(1 -q_{12}q_{21}q_{31}q_{13}) q_{21}q_{22}\}q_{12}^{ -1}q_{13}^{ -1}  [1, 3]  [1, 2]2 $
 
$+q_{13}^{ -1}q_{32}(1 -q_{12}q_{21}q_{31}q_{13} -q_{12}q_{21}q_{22}q_{31}q_{13}  + q_{12}^{2}q_{21}^{2}q_{22}q_{31}q_{13}) 2  [1, 3]  [1, 2]$
 
$ +q_{13}^{ -1}q_{32}q_{31}\{(q_{12}^{ -1} -q_{21}) -(1 -q_{12}q_{21}q_{31}q_{13}) q_{21}q_{22}\}  [1, 2]  [1, 3]2$
  
$+q_{13}^{ -1}q_{21}q_{22}q_{32}^{2}q_{12}q_{31}(1 -q_{12}q_{21}q_{31}q_{13} -q_{12}q_{21}q_{22}q_{31}q_{13}  
+q_{12}^{2}q_{21}^{2}q_{22}q_{31}q_{13})   [1, 2]2 [1, 3] $.

\end {Lemma}


According to  \cite [Table 2]{He05c}, the first node, second node and third node of every generalized Dynkin diagram 
denote $q_{33}, q_{11}, q_{22}$, respectively. Let $\mathbb B_V$ be the set of all hard super-letters in $\mathfrak B(V) $
(i.e. the generators of PBW basis.  Hard super-letters were defined in \cite [Def. 6] {{Kh99b}}) .

\begin{Theorem}\label{2.2'}

{\rm (i)}  If   $\begin{picture}(100,      20)  \put(27,     1) {\makebox(0,
0) [t]{$\bullet$}} \put(60,      1) {\makebox(0,     0) [t]{$\bullet$}}
\put(93,      1) {\makebox(0,      0) [t]{$\bullet$}} \put(28,    -1) {\line(1,
0) {33}} \put(61,     -1) {\line(1,      0) {30}} \put(18,      7) {$-1$} \put(35,      6) {$q$} \put(58,
7) {$-1$}  \put(75,      6) {$q^{-1}$}  \put(93,      7) {$-1$} \put(110,      1) {$,   $}\ \ \ \put(120,      1)  {$ q \in R_m,   m>2,$}
\end{picture}$ \\
then  $\mathbb B_V =\{   [x_1],     [x_2],     [x_3],     [x_1,   x_2],     [x_1,   x_3],
  [  [x_1,   x_3],   x_2]\}$  and $\dim \mathfrak B(V) =2^{4}m^{2}$.

{\rm (ii)}  If $\begin{picture}(100,      20)  \put(27,     1) {\makebox(0,
0) [t]{$\bullet$}} \put(60,      1) {\makebox(0,     0) [t]{$\bullet$}}
\put(93,      1) {\makebox(0,      0) [t]{$\bullet$}} \put(28,    -1) {\line(1,
0) {33}} \put(61,     -1) {\line(1,      0) {30}} \put(18,      7) {$-1$} \put(35,      6) {$\zeta $} \put(58,
7) {$-1$}  \put(75,      6) {$\zeta$}  \put(93,      7) {${-1}$} \put(110,      1) {$,   $}\ \ \ \put(120,      1)  {$ \zeta \in R_3, $}
\end{picture}$ \\
 then 
$\mathbb B_V=\{   [x_1],     [x_2],     [x_3],     [x_1,   x_2],     [x_1,   x_3],
  [  [x_1,   x_3],   x_2], [  [x_1,   x_2],    [x_1,   x_3]], $
 
$  [  [x_1,   x_2],     [  [x_1,   x_3],   x_2]],  [  [x_1,   x_3],    [  [x_1,   x_3],   x_2]],    
[  [  [x_1,   x_2],    [  [x_1,   x_3],   x_2]],    [x_1,   x_3]]\}$
\\ and    $\dim \mathfrak B(V) =2^{7}3^{4}$.

{\rm (iii)}   If  $\begin{picture}(100,      20)  \put(27,     1) {\makebox(0,
0) [t]{$\bullet$}} \put(60,      1) {\makebox(0,     0) [t]{$\bullet$}}
\put(93,      1) {\makebox(0,      0) [t]{$\bullet$}} \put(28,    -1) {\line(1,
0) {33}} \put(61,     -1) {\line(1,      0) {30}} \put(18,      7) {q} \put(35,      6) {$q^{-1}$} \put(58,   7) {$-1$}  \put(75,     
 6) {$r^{-1}$}  \put(93,      7) {$r$} \put(110,      1) {$,   $}\ \ \ \put(120,     
 1) {$ q \in R_m,      r \in R_{m'},   q\not= r, rq\not= 1; m,   m'>1$,} \end{picture}$ \\
then $\mathbb B_V=\{   [ x_1],     [x_2],     [x_3],     [x_1,   x_2],     [x_1,   x_3],
  [  [x_1,   x_3],   x_2],
 [  [x_1,   x_2],    [x_1,   x_3]]\}$ 
 and \\ $\dim \mathfrak B(V)  = 2^{4} \frac {m^2 m'{}^2}{(m,   m') }$.

\end {Theorem}

\begin{Proof}
  Assume that $  [u]$ is a hard super-letter or zero and  $u=vw$ is the Shirshow decomposition of $u$ when $[u] \not=0$. 
Applying Lemma \ref {2.2'''} we can  show  $  [u] \in \mathbb B_V$ step by step for the length $\mid u \mid $ of $u$. 
\end{Proof}

\begin{Lemma}\label{2.1} Let $n=km$ and $(s,  m)=1.$ $t_1,  t_2,  t_3\in \mathbb Z$.
 If  $t_1 \equiv 1 \ \ \ (\mbox {mod } n) $,  then the following conditions are equivalent.

 {\rm (i)} \begin {eqnarray}\label {e2.1.11'} \left \{
\begin
{array} {lll}  x_1 y_1 & \equiv  & t_1 sk \ \ \ (\mbox {mod } n)\\
 x_2 y_2 & \equiv  & t_2 sk \ \ \ (\mbox {mod } n)\\
 x_3 y_3 & \equiv  & t_3 sk \ \ \ (\mbox {mod } n)\\
 x_1 y_2+ x_2y_1 & \equiv  & t_4 sk \ \ \ (\mbox {mod } n)\\
 x_1 y_3+ x_3y_1 & \equiv  & t_5 sk \ \ \ (\mbox {mod } n)\\
x_3 y_2+ x_2y_3 & \equiv  & t_6 sk \ \ \ (\mbox {mod } n)
\end {array} \right.  \end {eqnarray} has a solution

{\rm (ii)}
 \begin {eqnarray}\label {e2.1.12} \left \{
\begin
{array} {llll}
 t_1(x_2)^2-t_4x_2 +t_2&\equiv &0  &\ \ \ (\mbox {mod } m)\\
 t_1(x_3)^2-t_5x_3 +t_3 & \equiv &0 &\ \ \ (\mbox {mod } m)\\
x _1 &\equiv&  1 &\ \ \ (\mbox {mod } n)\\
y_1&\equiv &t_1ks & \ \ \ (\mbox {mod } n)\\
y_2 & \equiv  &(t_4-x_2t_1)ks &\ \ \ (\mbox {mod } n)\\
 y_3 & \equiv  &(t_5-x_3t_1)ks &\ \ \ (\mbox {mod } n)\\
 x_1 y_1 & \equiv  & t_1 sk &\ \ \ (\mbox {mod } n)\\
 x_2 y_2 & \equiv  & t_2 sk &\ \ \ (\mbox {mod } n)\\
 x_3 y_3 & \equiv  & t_3 sk &\ \ \ (\mbox {mod } n)\\
 x_1 y_2+ x_2y_1 & \equiv  & t_4 sk &\ \ \ (\mbox {mod } n)\\
x_1 y_3+ x_3y_1 & \equiv  & t_5 sk &\ \ \ (\mbox {mod } n)\\
x_3 y_2+ x_2y_3 & \equiv  & t_6 sk &\ \ \ (\mbox {mod } n)\\
\end {array} \right.  \end {eqnarray} has a solution.

{\rm (iii)} \begin {eqnarray}\label {e2.1.51} \left \{
\begin
{array} {llll}
 t_1(x_2)^2-t_4x_2 +t_2&\equiv &0  &\ \ \ (\mbox {mod } m)\\
 t_1(x_3)^2-t_5x_3 +t_3 & \equiv &0 &\ \ \ (\mbox {mod } m)\\
x _1 &\equiv&  1 &\ \ \ (\mbox {mod } n)\\
y_1&\equiv &t_1ks & \ \ \ (\mbox {mod } n)\\
y_2 & \equiv  &(t_4-x_2t_1)ks &\ \ \ (\mbox {mod } n)\\
 y_3 & \equiv  &(t_5-x_3t_1)ks &\ \ \ (\mbox {mod } n)\\
2t_1x_2 x_3- t_4x_3-t_5 x_2  & \equiv  &-t_6  &\ \ \ (\mbox {mod } m)\\
\end {array} \right.  \end {eqnarray} has a solution.

\end{Lemma}

\begin{Lemma}\label{2.1''} Let $n=km$ and $(s,  m)=1;$ $t_1,  t_2,  t_3\in \mathbb Z$.
  If $(m, t_1) =1$, then   \begin {eqnarray}\label {e2.1.15} \left \{
\begin
{array} {lll}  x_1 y_1 & \equiv  & t_1 sk \ \ \ (\mbox {mod } n)\\
 x_2 y_2 & \equiv  & t_2 sk \ \ \ (\mbox {mod } n)\\
 x_3 y_3 & \equiv  & t_3 sk \ \ \ (\mbox {mod } n)\\
 x_1 y_2+ x_2y_1 & \equiv  & t_4 sk \ \ \ (\mbox {mod } n)\\
 x_1 y_3+ x_3y_1 & \equiv  & t_5 sk \ \ \ (\mbox {mod } n)\\
x_3 y_2+ x_2y_3 & \equiv  & t_6 sk \ \ \ (\mbox {mod } n)
\end {array} \right.  \end {eqnarray} has a solution if and only if
 \begin {eqnarray}\label {e2.1.31} \left \{
\begin
{array} {llll}
 t_1(x_2)^2-t_4x_2 +t_2&\equiv &0  &\ \ \ (\mbox {mod } m)\\
 t_1(x_3)^2-t_5x_3 +t_3 & \equiv &0 &\ \ \ (\mbox {mod } m)\\
x _1 &\equiv&  1 &\ \ \ (\mbox {mod } m)\\
y_1&\equiv &t_1& \ \ \ (\mbox {mod } m)\\
y_2 & \equiv  &(t_4-x_2t_1) &\ \ \ (\mbox {mod } m)\\
 y_3 & \equiv  &(t_5-x_3t_1)&\ \ \ (\mbox {mod } m)\\
 x_1 y_2+ x_2y_1 & \equiv  & t_4  &\ \ \ (\mbox {mod } m)\\
x_1 y_3+ x_3y_1 & \equiv  & t_5  &\ \ \ (\mbox {mod } m)\\
x_3 y_2+ x_2y_3 & \equiv  & t_6  &\ \ \ (\mbox {mod } m)\\
\end {array} \right.  \end {eqnarray} has a solution.

 \end{Lemma}

\begin{Lemma}\label{2.2''} Let $f$ denote the lowest common multiple  of $m$ and $m'$ with $(s, m)=1 = (s', m') $ and $m, m'>1.$ Then
 \begin {eqnarray}  \label {e2.4.1}\left \{
\begin
{array} {lll}  x_1 y_1 & \equiv  &  \frac {n}{2} \ \ \ (\mbox {mod } n)\\
 x_2 y_2 & \equiv  &   \frac {sn}{m} \ \ \ (\mbox {mod } n)\\
 x_3 y_3 & \equiv  &   \frac {s'n}{m'} \ \ \ (\mbox {mod } n)\\
 x_1 y_2+ x_2y_1 & \equiv  &-  \frac {sn}{m} \ \ \ (\mbox {mod } n)\\
 x_1 y_3+ x_3y_1 & \equiv  & -  \frac {s'n}{m'} \ \ \ (\mbox {mod } n)\\
x_3 y_2+ x_2y_3 & \equiv  & 0 \ \ \ (\mbox {mod } n)
\end {array} \right.  \end {eqnarray} has a solution
if and only if $\alpha _i = \alpha _i'$ when $\alpha _i \alpha _i' \not=0$ for $1\le i \le t$, and
\begin {eqnarray} \label {e2.4.70}
  -s   \equiv   m''s' \ \ \ (\mbox {mod } m') \end {eqnarray} when $m= m''m'$ and $(m', m'')=1$;
\begin {eqnarray} \label {e2.4.70'}
  -s '  \equiv   m''s \ \ \ (\mbox {mod } m) \end {eqnarray} when $m'= m''m$ and  $(m, m'')=1$. 
Here $m = 2^{\alpha_1}3^{\alpha_2} p_3 ^{\alpha_3}\cdots p_t^{\alpha_t}$,
$m '= 2^{\alpha_1'}3^{\alpha_2'} p_3 ^{\alpha_3'}\cdots p_t^{\alpha_t'}$ be the prime
decomposition, respectively.
\end {Lemma}

\begin{Theorem}\label{2.2}  If $(V, (q_{ij})_{3\times 3})$ is a braided vector space, then $V$ is a connected $\mathbb Z_n$-{\rm YD} 
module such that $\dim \mathfrak B(V)< \infty$ if and only if  one of  the following conditions holds:

({\rm i})   The generalized Dynkin diagram of $V$  is Weyl equivalent to \\
$\begin{picture}(100,    20) \put(27,   1){\makebox(0,
0)[t]{$\bullet$}} \put(60,    1){\makebox(0,   0)[t]{$\bullet$}}
\put(93,    1){\makebox(0,    0)[t]{$\bullet$}} \put(28,  -1){\line(1,
0){33}} \put(61,   -1){\line(1,    0){30}} \put(18,    7){$-1$} \put(35,    6){$q$} \put(58,
7){$-1$}  \put(75,    6){$q^{-1}$}  \put(93,    7){$-1$} \put(110,    1){$, $}\ \ \ \put(120,    1) {$ q \in R_m, m>2.$}
\end{picture}$ \\


({\rm ii})    The generalized Dynkin diagram of $V$  is Weyl equivalent to \\ $\begin{picture}(100,    20) \put(27,   1){\makebox(0,
0)[t]{$\bullet$}} \put(60,    1){\makebox(0,   0)[t]{$\bullet$}}
\put(93,    1){\makebox(0,    0)[t]{$\bullet$}} \put(28,  -1){\line(1,
0){33}} \put(61,   -1){\line(1,    0){30}} \put(18,    7){$-1$} \put(35,    6){$\zeta $} \put(58,
7){$-1$}  \put(75,    6){$\zeta$}  \put(93,    7){${-1}$} \put(110,    1){$, $}\ \ \ \put(120,    1) {$ \zeta \in R_3.$}
\end{picture}$ \\

 ({\rm iii})  The generalized Dynkin diagram of $V$  is Weyl equivalent to  \\ $\begin{picture}(100,    20) \put(27,   1){\makebox(0,
0)[t]{$\bullet$}} \put(60,    1){\makebox(0,   0)[t]{$\bullet$}}
\put(93,    1){\makebox(0,    0)[t]{$\bullet$}} \put(28,  -1){\line(1,
0){33}} \put(61,   -1){\line(1,    0){30}} \put(18,    7){q} \put(35,    6){$q^{-1}$} \put(58, 7){$-1$}  \put(75,   
 6){$r^{-1}$}  \put(93,    7){$r$} \put(110,    1){$, $}\ \ \ \put(120,    1) \\ \
\end{picture}$;  $\alpha _i = \alpha _i'$ when $\alpha _i \alpha _i' \not=0$ for $1\le i \le t;$ $
  -s   \equiv   m''s' \ \ \ (\mbox {mod } m') $ when $m= m''m'$ and $(m'', m')=1$;
$
  -s   \equiv   m''s' \ \ \ (\mbox {mod } m)$ when $m'= m''m$ and $(m, m'')=1$;
Here
$   q \in R_m, r\in R_{m'},   \omega  \in R_{n},$   $ m>1, m'>1; q\not= r, q\not= r^{-1};$ $(s, m)=1$; $(s', m')=1$; 
$q = \omega ^{\frac {ns}{m}}$, $r = \omega ^{\frac {ns'}{m'}}$;  $m = 2^{\alpha_1}3^{\alpha_2} p_3 ^{\alpha_3}\cdots p_t^{\alpha_t}$,
$m '= 2^{\alpha_1'}3^{\alpha_2'} p_3 ^{\alpha_3'}\cdots p_t^{\alpha_t'}$ be the prime decomposition, respectively.

\end{Theorem}

\begin{Proof} {\em The necessity}. By \cite [Th. 12]{He05c}, 
we  only need  to consider the generalized Dynkin diagrams in \cite [Table 2]{He05c}. 
The Dynkin diagrams above are in Row 8, 9,  15 of \cite [Table 2]{He05c}. So we need to exclude the Dynkin diagrams
in all other Rows of \cite [Table 2]{He05c}. This follows from the application of Lemma \ref {1.0} and Lemma \ref {2.1}.
For instance,  Row 1 of  \cite [Table 2]{He05c}. By Lemma \ref {1.0},
 \begin {eqnarray*} \left \{
\begin
{array} {lll}  x_1 y_1 & \equiv  &  sk \ \ \ (\mbox {mod } n)\\
 x_2 y_2 & \equiv  &  sk \ \ \ (\mbox {mod } n)\\
 x_3 y_3 & \equiv  &  sk \ \ \ (\mbox {mod } n)\\
 x_1 y_2+ x_2y_1 & \equiv  &- sk \ \ \ (\mbox {mod } n)\\
 x_1 y_3+ x_3y_1 & \equiv  & - sk \ \ \ (\mbox {mod } n)\\
x_3 y_2+ x_2y_3 & \equiv  & 0 \ \ \ (\mbox {mod } n)
\end {array} \right.  \end {eqnarray*} has a solution. Thus by Lemma \ref {2.1}
 \begin {eqnarray*} \left \{
\begin
{array} {llll}
x_1 &\equiv& 1 &\ \ \ (\mbox {mod } n)\\
y_1&\equiv &ks & \ \ \ (\mbox {mod } n)\\
y_2 & \equiv  &(-1- x_2)ks &\ \ \ (\mbox {mod } n)\\
 y_3 & \equiv  &(-1-x_3)ks &\ \ \ (\mbox {mod } n)\\
(x_2)^2 +x_2 +1&\equiv &0  &\ \ \ (\mbox {mod } m)\\
 (x_3)^2 +x_3 +1 & \equiv &0 &\ \ \ (\mbox {mod } m)\\
  2x_2x_3 + x_2+x_3 & \equiv &0 &\ \ \ (\mbox {mod } m)\\
 \end {array} \right.   \end {eqnarray*} has a solution,
 which implies that $2 \nmid m$ and \begin {eqnarray*} \left \{
\begin
{array} {llll}
(2x_2+1)^2&\equiv &-3  &\ \ \ (\mbox {mod } m)\\
 (2x_3+1)^2&\equiv &-3  &\ \ \ (\mbox {mod } m)\\
 (2x_2+1)(2x_3+1) &\equiv &1 &\ \ \ (\mbox {mod } m)\\
 \end {array} \right.  \end {eqnarray*} has a solution.  One gets  $9 \equiv 1 \ \ \ (\mbox {mod } m)$,  
which is a contradiction. So the diagram in Row 1 of \cite [Table 2]{He05c} is excluded.
By similar procedure, we can exclude the generalized Dynkin diagrams in all other Rows 
except those in Rows 8, 9, 15 of \cite [Table 2]{He05c} .

{\em The sufficiency}. It follows from Lemma \ref {2.2'} that $\dim \mathfrak B (V) < \infty$ 
when the generalized Dynkin diagrams are in Row 8, 9,  15 of \cite [Table 2]{He05c}.  By \cite [Th. 12]{He05c},  
we   need  to decide if   Row 8,  Row 9  and Row 15   in \cite [Table 2]{He05c} are $kG$- {\rm YD} modules.

({\rm i})
 Row 8 \cite [Table 2]{He05c}. There exists  a DDG  \begin{picture}(100,   20) \put(27,  1){\makebox(0,
0)[t]{$\bullet$}} \put(60,   1){\makebox(0,  0)[t]{$\bullet$}}
\put(93,   1){\makebox(0,   0)[t]{$\bullet$}} \put(28, -1){\line(1,
0){33}} \put(61,  -1){\line(1,   0){30}} \put(18,   7){q} \put(35,   6){$q^{-1}$} \put(58,
7){$-1$}  \put(75,   6){$q$}  \put(93,   7){$q^{-1}$} \put(110,   1){$,$}\ \ \ \put(120,   1) {} 
\end{picture}  \ \ \ \ \ \ \ \ \ \ \ \ \ \ \ \ \ \
{$ q \in R_m,$}  in Row 8 \cite [Table 2]{He05c}.  It follows from Lemma \ref {2.2''} when one sets $s=-s'$.

({\rm ii})  Row 15 \cite [Table 2]{He05c}. There exists  a DDG  \begin{picture}(100,   20) \put(27,  1){\makebox(0,
0)[t]{$\bullet$}} \put(60,   1){\makebox(0,  0)[t]{$\bullet$}}
\put(93,   1){\makebox(0,   0)[t]{$\bullet$}} \put(28, -1){\line(1,
0){33}} \put(61,  -1){\line(1,   0){30}} \put(18,   7){-1} \put(35,   6){$\xi^{-1}$} \put(58,
7){$\xi$}  \put(75,   6){$\xi$}  \put(93,   7){$-1$} \put(110,   1){$,$}\ \ \ \put(120,   1) {} 
\end{picture}  \ \ \ \ \ \ \ \ \ \ \ \ \ \ \ \ \ \
{$ \xi \in R_3,$}  in Row 15 \cite [Table 2]{He05c}.

 \begin {eqnarray*} \left \{
\begin{array} {lll}  x_1 y_1 & \equiv  & 2 sk \ \ \ (\mbox {mod } 6k)\\
 x_2 y_2 & \equiv  & 3sk \ \ \ (\mbox {mod } 6k)\\
 x_3 y_3 & \equiv  &  3 sk \ \ \ (\mbox {mod } 6k)\\
 x_1 y_2+ x_2 y_1 & \equiv  & 2 sk \ \ \ (\mbox {mod } 6k)\\
x_3 y_1+ x_1y_3 & \equiv  &-2sk \ \ \ (\mbox {mod } 6k)\\
x_2 y_3+ x_3 y_2 & \equiv  &  0 \ \ \ (\mbox {mod } 6k)
\end {array} \right.  
\end {eqnarray*}
has a solution: $x_2 = 1,  y_2 = 3ks, $ $ x_1 = 4,  y_1 =  2s{k}$,  $x_3 = 5,  y_3= 3{k}s$.

({\rm iii}) Row 9 \cite [Table 2]{He05c}. It follows from Lemma \ref {2.2''}.
\end{Proof}

\section{ Nichols algebras of diagonal type with rank $>3 $}\label {s3}

In this section we prove that finite dimensional  Nichols
algebra over $\mathbb Z_2$ is a quantum linear space and Nichols algebra of
connected  $\mathbb Z_n$-{\rm YD} module $V$  with $\dim V >3$ is infinite dimensional.

\begin{Theorem}\label{3.1}  If $V$ is a connected $k\mathbb Z_n$-Yetter-Drinfeld module with diagonal type and $rank >3$,
then  $ \dim \mathfrak B (V) = \infty$ and $\Delta (\mathfrak B (V))$ is infinite.
\end{Theorem}

\begin{Proof} It is enough to show this is the case for $\dim V =4$.

Except Row 18,  Row 20,  Row 21,  Row 22,
all GDDs  in  \cite [  Table B] {He09}   contain GDDs in  \cite [Table 2]{He05c}.
By Theorem \ref {2.2}, these four cases are not GDDs of  any $kG$-{\rm YD} modules.

({\rm i}) Row 18. $n = mk $,  $m =6$,  $(s,  m)=1.$ 
By Lemma \ref {1.0},
 \begin {eqnarray*} \left \{
\begin
{array} {lll}  x_1 y_1 & \equiv  & -2sk \ \ \ (\mbox {mod } n)\\
 x_2 y_2 & \equiv  & -2sk \ \ \ (\mbox {mod } n)\\
 x_3 y_3 & \equiv  &  2sk \ \ \ (\mbox {mod } n)\\
 x_1 y_2+ x_2y_1 & \equiv  & 2 sk \ \ \ (\mbox {mod } n)\\
 x_1 y_3+ x_3y_1 & \equiv  &  0 \ \ \ (\mbox {mod } n)\\
x_3 y_2+ x_2y_3 & \equiv  & 2sk \ \ \ (\mbox {mod } n)
\end {array} \right.  \end {eqnarray*} has a solution.  Let $s_1=2s$. Obviously,  $(s_1, 3)=1$. Thus \begin {eqnarray*} \left \{
\begin{array} {lll}  x_1 y_1 & \equiv  & -s_1k \ \ \ (\mbox {mod } 3k)\\
 x_2 y_2 & \equiv  & -s_1k \ \ \ (\mbox {mod } 3k)\\
 x_3 y_3 & \equiv  &  s_1k \ \ \ (\mbox {mod } 3k)\\
 x_1 y_2+ x_2y_1 & \equiv  & s_1k \ \ \ (\mbox {mod } 3k)\\
 x_1 y_3+ x_3y_1 & \equiv  &  0 \ \ \ (\mbox {mod } 3k)\\
x_3 y_2+ x_2y_3 & \equiv  & s_1k \ \ \ (\mbox {mod } 3k)
\end {array} \right.  \end {eqnarray*} has a solution.
Thus by Lemma \ref {2.1}
 \begin {eqnarray*} \left \{
\begin{array} {llll}
x_3 &\equiv& 1 &\ \ \ (\mbox {mod } 3k)\\
y_3&\equiv & s_1k  & \ \ \ (\mbox {mod } 3k)\\
y_1 & \equiv  &- x_1 s_1k  &\ \ \ (\mbox {mod } 3k)\\
 y_2 & \equiv  &(1-x_2) s_1k  &\ \ \ (\mbox {mod } 3k)\\
(x_1)^2-1 &\equiv &0  &\ \ \ (\mbox {mod } 3)\\
 (x_2)^2-x_2-1 & \equiv &0 &\ \ \ (\mbox {mod } 3)\\
 -2x_1x_2 +x_1 & \equiv &1 &\ \ \ (\mbox {mod } 3)\\
 \end {array} \right.   \end {eqnarray*}has a solution,
 which implies that  \begin {eqnarray*} \left \{
\begin{array} {llll}
(2x_2-1)^2&\equiv &5 &\ \ \ (\mbox {mod } 3)\\
 (x_1)^2&\equiv &1  &\ \ \ (\mbox {mod } 3)\\
 x_1(-2x_2+1) &\equiv &1 &\ \ \ (\mbox {mod } 3)\\
 \end {array} \right..  \end {eqnarray*}  One gets $5 \equiv 1 \ \ \ (\mbox {mod } 3)$,   which is a contradiction.

({\rm ii}) Row 20.
$n = mk $,  $m =6$,  $(s,  m)=1.$ Consider the last GDD in Row 21.
By Lemma \ref {1.0},
 \begin {eqnarray*} \left \{
\begin{array} {lll}  x_1 y_1 & \equiv  &  2sk \ \ \ (\mbox {mod } n)\\
 x_2 y_2 & \equiv  & 2sk \ \ \ (\mbox {mod } n)\\
 x_3 y_3 & \equiv  & -2sk \ \ \ (\mbox {mod } n)\\
 x_1 y_2+ x_2y_1 & \equiv  &-2 sk \ \ \ (\mbox {mod } n)\\
 x_1 y_3+ x_3y_1 & \equiv  &  0 \ \ \ (\mbox {mod } n)\\
x_3 y_2+ x_2y_3 & \equiv  & 2sk \ \ \ (\mbox {mod } n)
\end {array} \right.  \end {eqnarray*} has a solution.  Let $s_1=2s$. Obviously,  $(s_1, 3)=1$. Thus \begin {eqnarray*} \left \{
\begin{array} {lll}  x_1 y_1 & \equiv  &  s_1k \ \ \ (\mbox {mod } 3k)\\
 x_2 y_2 & \equiv  &  s_1k \ \ \ (\mbox {mod } 3k)\\
 x_3 y_3 & \equiv  &- s_1k \ \ \ (\mbox {mod } 3k)\\
 x_1 y_2+ x_2y_1 & \equiv  &-s_1k \ \ \ (\mbox {mod } 3k)\\
 x_1 y_3+ x_3y_1 & \equiv  &  0 \ \ \ (\mbox {mod } 3k)\\
x_3 y_2+ x_2y_3 & \equiv  & s_1k \ \ \ (\mbox {mod } 3k)
\end {array} \right.  \end {eqnarray*} has a solution.
Thus by Lemma \ref {2.1}
 \begin {eqnarray*} \left \{
\begin{array} {llll}
x_1 &\equiv& 1 &\ \ \ (\mbox {mod } 3k)\\
y_1&\equiv & s_1k  & \ \ \ (\mbox {mod } 3k)\\
y_2 & \equiv  &(-1- x_2) s_1k  &\ \ \ (\mbox {mod } 3k)\\
 y_3 & \equiv  &(-x_3) s_1k  &\ \ \ (\mbox {mod } 3k)\\
(x_2)^2 +x_2+1 &\equiv &0  &\ \ \ (\mbox {mod } 3)\\
 (x_3)^2-1 & \equiv &0 &\ \ \ (\mbox {mod } 3)\\
  2x_2x_3 +x_3 & \equiv &-1 &\ \ \ (\mbox {mod } 3)\\
 \end {array} \right.   \end {eqnarray*}has a solution,
 which implies that  \begin {eqnarray*} \left \{
\begin
{array} {llll}
(2x_2+1)^2&\equiv &-3 &\ \ \ (\mbox {mod } 3)\\
 (x_3)^2&\equiv &1  &\ \ \ (\mbox {mod } 3)\\
 x_3(2x_2+1) &\equiv &-1 &\ \ \ (\mbox {mod } 3)\\
 \end {array} \right..  \end {eqnarray*}  One gets $-3 \equiv 1 \ \ \ (\mbox {mod } 3)$,   which is a contradiction.

({\rm iii})  Row 21. $n = mk $,  $m =6$,  $(s,  m)=1.$  Consider the last GDD in Row 21.
By Lemma \ref {1.0},
 \begin {eqnarray*} \left \{
\begin
{array} {lll}  x_1 y_1 & \equiv  &  2sk \ \ \ (\mbox {mod } n)\\
 x_2 y_2 & \equiv  & 2sk \ \ \ (\mbox {mod } n)\\
 x_3 y_3 & \equiv  &  2sk \ \ \ (\mbox {mod } n)\\
 x_1 y_2+ x_2y_1 & \equiv  &-2 sk \ \ \ (\mbox {mod } n)\\
 x_1 y_3+ x_3y_1 & \equiv  &  0 \ \ \ (\mbox {mod } n)\\
x_3 y_2+ x_2y_3 & \equiv  &-2sk \ \ \ (\mbox {mod } n)
\end {array} \right.  \end {eqnarray*} has a solution.  Let $s_1=2s$. Obviously,  $(s_1, 3)=1$. Thus \begin {eqnarray*} \left \{
\begin
{array} {lll}  x_1 y_1 & \equiv  &  s_1k \ \ \ (\mbox {mod } 3k)\\
 x_2 y_2 & \equiv  &  s_1k \ \ \ (\mbox {mod } 3k)\\
 x_3 y_3 & \equiv  &  s_1k \ \ \ (\mbox {mod } 3k)\\
 x_1 y_2+ x_2y_1 & \equiv  &-s_1k \ \ \ (\mbox {mod } 3k)\\
 x_1 y_3+ x_3y_1 & \equiv  &  0 \ \ \ (\mbox {mod } 3k)\\
x_3 y_2+ x_2y_3 & \equiv  &-s_1k \ \ \ (\mbox {mod } 3k)
\end {array} \right.  \end {eqnarray*} has a solution.
Thus by Lemma \ref {2.1}
 \begin {eqnarray*} \left \{
\begin
{array} {llll}
x_1 &\equiv& 1 &\ \ \ (\mbox {mod } 3k)\\
y_1&\equiv & s_1k  & \ \ \ (\mbox {mod } 3k)\\
y_2 & \equiv  &(-1- x_2) s_1k  &\ \ \ (\mbox {mod } 3k)\\
 y_3 & \equiv  &(-x_3) s_1k  &\ \ \ (\mbox {mod } 3k)\\
(x_2)^2 +x_2+1 &\equiv &0  &\ \ \ (\mbox {mod } 3)\\
 (x_3)^2 +1 & \equiv &0 &\ \ \ (\mbox {mod } 3)\\
  2x_2x_3 +x_3 & \equiv &1 &\ \ \ (\mbox {mod } 3)\\
 \end {array} \right.   \end {eqnarray*}has a solution,
 which implies that  \begin {eqnarray*} \left \{
\begin
{array} {llll}
(2x_2+1)^2&\equiv &-3 &\ \ \ (\mbox {mod } 3)\\
(x_3)^2 +1 & \equiv &0 &\ \ \ (\mbox {mod } 3)\\
 x_3(2x_2+1) &\equiv &1 &\ \ \ (\mbox {mod } 3)\\
 \end {array} \right..  \end {eqnarray*}  One gets $0 \equiv 1 \ \ \ (\mbox {mod } 3)$,   which is a contradiction.

({\rm iv}) Row 22. $n = mk $,  $m =4$,  $(s,  m)=1.$ By Lemma \ref {1.0},
 \begin {eqnarray*}\left \{
\begin
{array} {llll}  x_3 y_3 & \equiv  &  sk &\ \ \ (\mbox {mod } n)\\
 x_4 y_4& \equiv  & 3 sk &\ \ \ (\mbox {mod } n)\\
 x_4 y_3+ x_3y_4 & \equiv  & sk &\ \ \ (\mbox {mod } n)
\end {array} \right.  \end {eqnarray*}has a solution.
By Lemma \ref {1.2} ({\rm i}),  \begin {eqnarray*}\left .
\begin
{array} {lll} x^2-x +3 \equiv 0 & &\ \ \ (\mbox {mod } 4)\\
\end {array} \right. \end {eqnarray*}has a solution, which contradicts to Lemma \ref {4.2}({\rm i}).

\end{Proof}

\begin{Corollary}\label{3.2} ({\rm i}) If $V$ is a connected finite dimensional {\rm YD} module over $\mathbb Z_n$
 with $\dim \mathfrak B(V)< \infty$, then $\dim V <4.$

({\rm ii}) If $V$ is a finite dimensional {\rm YD} module over $\mathbb Z_n$ with $\dim \mathfrak B(V)< \infty$, 
then dimension of every connected component of $V$ is lesser than $4$.
\end{Corollary}

\begin{Proof} ({\rm i}) It follows from  Theorem \ref {3.1}.

({\rm ii}) It follows from Part ({\rm i}) and \cite [Lemma 4.2] {AS00}. 
\end{Proof}

\begin{Corollary}\label{3.3} If $V$ is a finite dimensional {\rm YD}  module over $\mathbb Z_n$ with braided matrix 
$(q_{ij})$ and ${\rm ord }(q_{ii})\not=1$,  then the following conditions are equivalent:

 ({\rm i}) $\dim \mathfrak B(V)< \infty$.

 ({\rm ii}) $\Delta (\mathfrak B(V))$ is finite.

 ({\rm iii}) $(\Delta (\mathfrak B(V)),  \chi_0,  E_0)$  is an arithmetic root system.

\end{Corollary}

The concept of  quantum linear spaces was introduced  in \cite [P673] {AS98}. In this case, $q_{ij}q_{ji}=1$ for $i \not= j.$
\begin{Corollary}\label{3.4} Every finite dimensional  Nichols algebra   over $\mathbb Z_2$ is a quantum linear space.
\end{Corollary}

\begin{Corollary}\label{3.5}  Assume that  $(V, (q_{ij})_{2\times 2})$ is a braided vector space.

{\rm (I)} If  $ p$ is  a prime number, then $V$ is a connected $\mathbb Z_p$-{\rm YD} module such that 
$\dim \mathfrak B(V)< \infty$  if and only if one of  the following conditions holds:

 T2(1) $1-q_{11}q_{12}q_{21} =1-q_{12}q_{21}q_{22}=0$,  $q_{12}q_{21}\in R_p$;
  $p=3$ or $p>3$ and   $(\frac {-3}{p})=1$.

 T2(2)$_1$ $1+q_{11} =1-q_{12}q_{21}q_{22}=0$,   $q_{12}q_{21}\in R_p$; $p>2.$

  T2(2)$_2$  $1+q_{22} =1-q_{12}q_{21}q_{11}=0$, $q_{12}q_{21}\in R_p$,
$p>2.$

  T2(3)  $1+q_{11} =1+q_{22}=0$,  $q_{12}q_{21}\in R_p$, $p>2.$

  T3(1)$_1$   $q_{12}q_{21} =q_{11}^{-2}$,  $q_{22} = q_{11}^2$,  $q_{11} \in R_p$; $p>3$ and  $p \equiv 1 \ \ \ (\mbox {mod } 4). $

 T3(1)$_2$ $q_{12}q_{21} =q_{11}^{-2}$,  $q_{22} =-1$, $ q_{11}\in {R_{p}}$;
 $p>2$.

T8(1) $ q_{12}q_{21}=q_{11}^{-3}$,  $
q_{22}=q_{11}^3$, $q_{11} \in R_{p}$, $p>3$ and  $(\frac {-3}{p})=1$.

{\rm (II)} Let $p$ be a prime number,  $n= p^\beta$ and $m = p^{\alpha}$  with  $0<\alpha \le \beta$ and $\beta >1$. 
Then $V$ is a connected $\mathbb Z_n$-{\rm YD} module such that $\dim \mathfrak B(V)< \infty$  
if and only if one of  the following conditions holds:

 T2(1) $1-q_{11}q_{12}q_{21} =1-q_{12}q_{21}q_{22}=0$,  $q_{12}q_{21}\in R_m$; $p=3, \alpha =1;$
  $p>3$ and   $(\frac {-3}{p})=1$.

 T2(2)$_1$ $1+q_{11} =1-q_{12}q_{21}q_{22}=0$,   $q_{12}q_{21}\in R_m$; $p=2, \alpha >1;$ $p>2.$

  T2(2)$_2$  $1+q_{22} =1-q_{12}q_{21}q_{11}=0$, $q_{12}q_{21}\in R_m$; $p=2, \alpha >1;$ $p>2.$

  T2(3)  $1+q_{11} =1+q_{22}=0$,  $q_{12}q_{21}\in R_m$; $p=2, \alpha >1;$ $p>2.$

  T3(1)$_1$   $q_{12}q_{21} =q_{11}^{-2}$,  $q_{22} = q_{11}^2$,  $q_{11} \in R_m$,  $m >2$; $p>3$ and  $p \equiv 1 \ \ \ (\mbox {mod } 4). $

 T3(1)$_2$ $q_{12}q_{21} =q_{11}^{-2}$,  $q_{22} =-1$, $ q_{11}\in {R_{m}}$, $m>2$; $p=2$,  $\alpha>3;$
 $p>2$.

 T3(2)$_1$  $\omega \in R_n$, $s =1, 2;$ $q_{11}= \omega ^{\frac {ns}{3}}$,
 $q_{22}= \omega ^{\frac {n}{m}}$, $q_{12}q_{21} q_{22}=1$, $m>3;$ $p=3$ and $\alpha >1.$

T8(1) $ q_{12}q_{21}=q_{11}^{-3}$,  $
q_{22}=q_{11}^3$, $q_{11} \in R_{m}$, $m>3$;  $p>3$  and  $ (\frac {-3}{p} )=1$.

 T8(2)$_{1}$ $ (q_{12}q_{21})^4=-1$, $q_{22}=-1$, $q_{12}q_{21}=-q_{11}$;  $m=8$;  $\alpha =3$.

  T8(2)$_{2}$ $ (q_{12}q_{21})^4=-1$, $q_{22}=-1$, $q_{11} =(q_{12}q_{21})^{-2}$;  $m = 8$,  $\alpha =3$.

 T8(3) $ (q_{12}q_{21})^4=-1$, $q_{11} =(q_{12}q_{21})^2$, $q_{22} =(q_{12}q_{21})^{-1}$; $m = 8$,  $\alpha =3$.

\end{Corollary}

\begin{Proof} It follows from Theorem \ref {1.7}. \end{Proof}

\begin{Corollary}\label{3.6}  Assume that  $(V, (q_{ij})_{3\times 3})$ is a braided vector space.

{\rm (I)} If  $ p$ is  a prime number, then $V$ is a connected $\mathbb Z_p$-{\rm YD} module such that 
$\dim \mathfrak B(V)< \infty$  if and only if one of  the following conditions holds:

({\rm i})   The generalized Dynkin diagram of $V$  is Weyl equivalent to \\
$\begin{picture}(100,    20) \put(27,   1){\makebox(0,
0)[t]{$\bullet$}} \put(60,    1){\makebox(0,   0)[t]{$\bullet$}}
\put(93,    1){\makebox(0,    0)[t]{$\bullet$}} \put(28,  -1){\line(1,
0){33}} \put(61,   -1){\line(1,    0){30}} \put(18,    7){$-1$} \put(35,    6){$q$} \put(58,
7){$-1$}  \put(75,    6){$q^{-1}$}  \put(93,    7){$-1$} \put(110,    1){$, $}\ \ \ \put(120,    1) {$ q \in R_p; p>2.$}
\end{picture}$ \\

 ({\rm ii})  The generalized Dynkin diagram of $V$  is Weyl equivalent to  \\ $\begin{picture}(100,    20) \put(27,   1){\makebox(0,
0)[t]{$\bullet$}} \put(60,    1){\makebox(0,   0)[t]{$\bullet$}}
\put(93,    1){\makebox(0,    0)[t]{$\bullet$}} \put(28,  -1){\line(1,
0){33}} \put(61,   -1){\line(1,    0){30}} \put(18,    7){q} \put(35,    6){$q^{-1}$} \put(58, 7){$-1$}  \put(75,   
 6){$r^{-1}$}  \put(93,    7){$r$} \put(110,    1) \ \ \ {$ {} $}\ \ \ \put(120,    1) \\ \
\end{picture}$; $
  -s   \equiv   s' \ \ \ (\mbox {mod } p)$. Here
$   q, r,  \omega \in R_p,$  $ q\not= r, q\not= r^{-1};$ $(s, p)=1$; $(s', p)=1$; $q = \omega ^{s }$, $r = \omega ^{s'}$.

{\rm (II)} Let $p$ be a prime number,  $n= p^\beta$ and $m = p^{\alpha}$  with  $0<\alpha \le \beta$ and $\beta >1$. 
Then $V$ is a connected $\mathbb Z_n$-{\rm YD} module such that $\dim \mathfrak B(V)< \infty$ 
 if and only if one of  the following conditions holds:

({\rm i})   The generalized Dynkin diagram of $V$  is Weyl equivalent to \\
$\begin{picture}(100,    20) \put(27,   1){\makebox(0,
0)[t]{$\bullet$}} \put(60,    1){\makebox(0,   0)[t]{$\bullet$}}
\put(93,    1){\makebox(0,    0)[t]{$\bullet$}} \put(28,  -1){\line(1,
0){33}} \put(61,   -1){\line(1,    0){30}} \put(18,    7){$-1$} \put(35,    6){$q$} \put(58,
7){$-1$}  \put(75,    6){$q^{-1}$}  \put(93,    7){$-1$} \put(110,    1){$, $}\ \ \ \put(120,    
1) {$ q \in R_m; m>2$; $ p =2$ and $\alpha >1$ or $ p>2.$}
\end{picture}$ \\

 ({\rm ii})  The generalized Dynkin diagram of $V$  is Weyl equivalent to  \\ \begin{picture}(100,    20) \put(27,   1){\makebox(0,
0)[t]{$\bullet$}} \put(60,    1){\makebox(0,   0)[t]{$\bullet$}}
\put(93,    1){\makebox(0,    0)[t]{$\bullet$}} \put(28,  -1){\line(1,
0){33}} \put(61,   -1){\line(1,    0){30}} \put(18,    7){q} \put(35,    6){$q^{-1}$} \put(58, 7){$-1$}  \put(75,   
 6){$r^{-1}$}  \put(93,    7){$r$} \put(110,    1) \ \ \ {$ {} $}\ \ \ \put(120,    1) \\ \
\end{picture}; $  -s   \equiv   s' \ \ \ (\mbox {mod } m)$; $m'= m>1;$ $p=2$ and  $\alpha >1$ or  $ p >2$.   Here
$   q \in R_m, r\in R_{m'},   \omega  \in R_{n},$  $ q\not= r, q\not= r^{-1};$ $(s, m)=1$; $(s', m')=1$; 
$q = \omega ^{\frac {ns}{m}}$, $r = \omega ^{\frac {ns'}{m'}}$.

\end{Corollary}

\begin{Proof} It follows from Theorem \ref {2.2}. \end{Proof}

\appendix
\section {Appendix }\label {s54}
In this section,  we recall some results on solutions of equation systems in $\mathbb Z_n$
\cite {Hu67} and braided vector spaces.

\subsection {Equation systems  in $\mathbb Z_n$ }

If prime $p \nmid a$ and  $x^2 \equiv a  \ \ \ (\mbox {mod } p)$ has a solution,  then $a$ is called a quadratic residue  of module $p$. 
Set \begin {eqnarray}\label {e4.0.1} (\frac {a}{p}) := \left \{
\begin
{array} {lll}  1  & \hbox { when } a \hbox { is  a quadratic residue  of module } p \\
-1  & \hbox { when } a \hbox { is  a quadratic   non-residue of module } p \\
\end {array} \right.. \end {eqnarray} This is called  Legendre sign.

\begin{Lemma}\label{4.1} Let $f(x) = a_n x^n + \cdots+a_1x +a_0 \in \mathbb Z [x]$ and  
$f'(x) := na_n x^{n-1} +(n-1)a_{n-1} x^{n-2} + \cdots +a_1$.

({\rm i}) If   $f(x)\equiv 0\ \ \ (\mbox {mod } p) $ and $f'(x)\equiv 0\ \ \
(\mbox {mod } p) $ has not any  common solution with  prime number $p$,  then $f(x)\equiv 0\ \ \ (\mbox {mod } p^k)$ 
has a solution if and only if   $f(x)\equiv 0\ \ \ (\mbox {mod } p) $ has a solution.

({\rm ii})  $ax +b\equiv 0\ \ \ (\mbox {mod } m) $ has a solution if and only if
 $ (a, m) \mid b $.
\end {Lemma}

\begin{Lemma}\label{4.2} Let  \begin
{eqnarray}\label {e4.2.11} f(x) := a x^2 + bx +c \equiv 0  \ \ \ (\mbox {mod } p^k),
\end {eqnarray} with prime $p$,    $p \nmid (a, b, c)$ and $k \in \mathbb N$.

({\rm i}) If $2 \nmid a$,  $2 \nmid b$,  then $2 \mid c$ if and only
if (\ref {e4.2.11}) has a solution when $p=2$.

({\rm ii})  If   $2 \nmid a$ and  $2 \mid b$,  then (\ref {e4.2.11}) is equivalent to $(ax+\frac{b}{2})^2\equiv\frac{b^2}{4}-ac \
\ \ (\mbox {mod } 2^k) $ when $p=2$.

({\rm iii}) If $p>2$,  $p \mid a$,  $p \nmid b$,  then  (\ref {e4.2.11})
always has a solution.

({\rm iv}) If $p>2$,  $p \nmid a$,   then (\ref {e4.2.11}) is equivalent to $(2ax+b)^2\equiv b^2-4ac \ \ \ (\mbox {mod } p^k) $. 
Furthermore (\ref {e4.2.11}) has a solution if and only if
 $(2ax+b)^2\equiv b^2-4ac \ \ \ (\mbox {mod } p) $ has a solution.

\end {Lemma}

\begin{Lemma}\label{4.3} Let
 \begin {eqnarray}\label {e4.1.11} x^2 \equiv a   \ \ \ (\mbox {mod } p^k)
\end {eqnarray} where prime $p \nmid a$,
$k \in \mathbb N$.

({\rm i}) If $p>2$,  then the number of solution of (\ref {e4.1.11}) is $1+(\frac{a}{p})$.

({\rm ii}) If $p=2$,  then

(1)  (\ref {e4.1.11}) has a solution when $k=1$.

(2) (\ref {e4.1.11})  has  two solutions  when   $k=2$ and  $a \equiv 1 $\ \ \ (\mbox {mod } $4) $.

(3) (\ref {e4.1.11})  has not any  solutions  when   $k=2$ and  $a \not\equiv 1 $\ \ \ (\mbox {mod } $4) $.

(4) (\ref {e4.1.11})  has  four  solutions  when  $k>2$ and  $a \equiv 1 \ \ \ (\mbox {mod } 8) $.

(5) (\ref {e4.1.11})  has not any   solutions  when  $k>2$ and  $a \not\equiv 1 \ \ \ (\mbox {mod } 8) $.

\end {Lemma}

\begin{Lemma}\label{4.4} Let $m=m_1m_2\cdots m_r$, where $
m_1,m_2, \cdots, m_r $ are pairwise relatives prime. then $f(x)\equiv 0\ \ \ (\mbox {mod } m ) $ has a solution
 if and only if every equation below has a solution : \begin {eqnarray*}\left.
\begin{array} {llll}
f(x) &\equiv & 0\ \ \ & (\mbox {mod } m_1) \\
f(x) &\equiv & 0\ \ \ & (\mbox {mod } m_2) \\
\cdots  &\cdots & \cdots\ \ \ &  \\
f(x) &\equiv & 0\ \ \ & (\mbox {mod } m_r)\\
\end {array} \right. . \end {eqnarray*}
\end {Lemma}

\subsection {Braided vector space}

If $\sigma \in \mathbb S_r$ and $q_{\sigma (i), \sigma (j)} = q_{ij}'$ for any $1\le i, j \le r,$ then the two matrixes 
$\left (
\begin{array} {llll} q_{11} & q_{12} & \cdots & q_{1r}\\
q_{21} & q_{22} & \cdots & q_{2r}\\
\cdots &\cdots &\cdots & \cdots\\
q_{r1} & q_{r2} & \cdots & q_{rr}\\
\end {array} \right) $  and $\left (
\begin
{array} {llll} q_{11} '& q_{12} '& \cdots & q_{1r}'\\
q_{21}' & q_{22} '& \cdots & q_{2r}'\\
\cdots &\cdots &\cdots & \cdots\\
q_{r1} '& q_{r2}' & \cdots & q_{rr}'\\
\end {array} \right) $  are called to be   permutation
similarity. In this case, GDDs of the two matrixes are called to be isomorphic.

If $(q_{ij})$ and $(q_{ij}')$ are permutation
similarity, then  the two braided vector spaces  $(V, (q_{ij}))$  and  $(V, (q_{ij}'))$ are the same since 
$x_{\sigma (1)}, x_{\sigma (2)}, \cdots, x_{\sigma (r)}$ is also a basis of $V$ with 
$C(x_{\sigma (i)} \otimes x_{\sigma (j)}) = q_{\sigma (i)\sigma (j)}  
x_{\sigma (j)} \otimes x_{\sigma (i)} = q'_{ij}x_{\sigma (j)} \otimes x_{\sigma (i)}$.

Recall \cite {ZZC04}.  $(G, \overrightarrow{g}, \overrightarrow{\chi}, J)$ is called an
{\it element system with characters} $($simply, {\rm ESC}$)$ if $G$
is a group, $J$ is a set, $\overrightarrow{ g } = \{g_i\} _{ i\in J}
\in Z(G)^J$ and $\overrightarrow{ \chi }= \{\chi_i\}_{ i \in J} \in
\widehat{G}^J $ with $ g_i \in Z(G)$ and $\chi _i \in \widehat G$.
${\rm ESC} (G, \overrightarrow{g}, \overrightarrow{\chi}, J)$ and
${\rm ESC} (G', \overrightarrow{g'}, \overrightarrow{\chi'}, J')$
are said to be isomorphic if there exist a group isomorphism $\phi:
G \rightarrow G'$ and a bijective map $\sigma: J\rightarrow J'$ such
that $\phi(g_i)=g'_{\sigma(i)}$ and $\chi'_{\sigma(i)}\phi=\chi_i$
for any $i \in J$.

Let $(G, g_i, \chi_i; i\in J)$ be an ${\rm ESC}$. Let $V$ be a
$k$-vector space with ${\rm dim}(V)=|J|$. Let $\{x_i \mid i\in J\}$
be a basis of $V$ over $k$. Define a left $kG$-action and a left
$kG$-coaction on $V$ by
$$g\cdot x_i = \chi_i(g)x_i,\ \delta^-(x_i )= g_i \otimes
x_i,\ i\in J,\ g\in G.$$ Then it is easy to see that $V$ is a
pointed {\rm YD} $kG$-module and $kx_i$ is a one dimensional {\rm
YD} $kG$-submodule of $V$ for any $i\in J$. Denote by $V(G, g_i,
\chi_i; i\in J)$ the pointed {\rm YD} $kG$-module $V$. Obviously,  $C(x_i\otimes x_j) = \chi _j(g_i) x_j \otimes x_i$ 
for any $i, j \in J,$ is the braiding.  (See \cite [Lemma 2.3 and Lemma 2.4]{ZZC04}) Every $kG$-{\rm YD} module is isomorphic to $V(G, g_i,
\chi_i; i\in J)$, which is a braided vector space with diagonal type and braided matrix $(q_{ij})= (\chi_j (g_i))$ when $J= \{1, 2, \cdots, r\}$.

\begin{Lemma}\label{4.5}

If  There is a Hopf
algebra isomorphism $\phi: kG\rightarrow kG'$ such that $V(G, g_i,
\chi_i;$ $i\in J)\cong\ ^{\phi^{-1}}_{\phi}V'(G' g_i', \chi_i';$
$i\in J')$ as {\rm YD} $kG$-modules with $J= J' = \{1, 2, \cdots, r\}$ and $G=G'$, then
  $(q_{ij})_{r\times r}$ and $(q_{ij}')_{r\times r}$ are permutation
similarity,
 where $q_{ij} = \chi _j (g_i)$ and $q_{ij}' = \chi _j '(g_i')$ for $1, 2, \cdots, r.$

 \end {Lemma}
 
\begin{Proof}   By \cite [Th. 4]{ZZC04},  ${\rm ESC}  (G, g_i, \chi_i; i\in J) \cong \ {\rm ESC}
(G', g_i', \chi_i'; i\in J')$ with $J= J' = \{1, 2, \cdots, r\}$. Consequently,
 there exists a bijective map $\sigma: J\rightarrow J'$ such
that $\phi(g_i)=g'_{\sigma(i)}$ and $\chi'_{\sigma(i)}\phi (g_j)=\chi_i (g_j)$
for any $i, j \in J$. That is, $q_{\sigma (j), \sigma (i)}' = q_{ji}$  for any $i, j \in J$. 
\end{Proof}

\begin{Corollary} \label{2.2''c}  Assume  $q_{11} =-1$ and   $(q_{22} +1) (q_{22} q_{12}q_{21}-1) =0.$  
If $V$ is connected  with rank 2, then the generators of PBW basis $\mathbb B_V = \{x_1,   x_2,   [x_1,   x_2]\}.$
\end {Corollary}

\begin{Proof} It follows  By Lemma \ref {2.2'''}.  \end{Proof}

\begin{Remark} In this paper,  the first node, second node and third node of every generalized Dynkin diagram 
denote $q_{33}, q_{11}, q_{22}$, respectively. For example,

\begin{picture}(100,    20)  \put(27,   1) {\makebox(0,
0) [t]{$\bullet$}} \put(60,    1) {\makebox(0,   0) [t]{$\bullet$}}
\put(93,    1) {\makebox(0,    0) [t]{$\bullet$}} \put(28,  -1) {\line(1,
0) {33}} \put(61,   -1) {\line(1,    0) {30}} \put(18,    7) {q} \put(35,    
6) {$q^{-1}$} \put(58, 7) {$-1$}  \put(75,    6) {$r^{-1}$}  \put(93,    
7) {$r$} \put(110,    1) {$, $}\ \ \ \put(120,    1) {$ q \in R_m,    r \in R_{m'}, q\not= r; m, m' >1$} 
\end{picture} \\  
denotes $q_{11} =-1, q_{22} =r, q_{33} =q, $ $q_{12}q_{21} = r^{-1}$, 
$q_{13}q_{31} = q^{-1}$.
\end{Remark}

\vskip 1.5cm
\noindent{\Large\bf Acknowledgement}
\vskip.1in
\noindent We would like to thank Prof. N. Andruskiewitsch for pointing out an error in the previous version,
which has led to the current substantially modified version.  Y.Z.Z was financially supported by the Australian Research Council.

\end {document}